\documentclass[11pt,a4paper]{amsart}

\usepackage{amssymb,amsmath,graphics,verbatim}
\usepackage{latexsym}
\usepackage{eucal}
\usepackage{a4wide}

\newtheorem{proposition}{Proposition}[section]
\newtheorem{theorem}{Theorem}[section]
\newtheorem{lemma}[theorem]{Lemma}

\newtheorem{coro}[theorem]{Corollary}
\newtheorem{remark}[theorem]{Remark}

\newcommand{\mc}{\mathcal}
\newcommand{\cal}{\mathcal}
\newcommand{\rr}{\mathbb{R}}
\newcommand{\R}{\mathbb{R}}
\newcommand{\C}{\mathbb{C}}
\newcommand{\nn}{\mathbb{N}}
\newcommand{\cc}{\mathbb{C}}

\newcommand{\zz}{\mathbb{Z}}

\newcommand{\eps}{\epsilon}

\newcommand{\pl}{\partial}
\newcommand{\x}{\times}

\newcommand{\til}{\widetilde}
\newcommand{\bbar}{\overline}

\newcommand{\cjd}{\rangle}
\newcommand{\cjg}{\langle}

\newcommand{\demi}{\frac{1}{2}}

\def\qed{\hfill$\square$\medskip}

\begin{document}
\title[Calder\'on Problem with partial data on Riemann surfaces]{Calder\'on inverse Problem with partial data on Riemann Surfaces}
\author{Colin Guillarmou}
\address{Laboratoire J.A. Dieudonn\'e\\
U.M.R. 6621 CNRS\\
Universit\'e de Nice Sophia-Antipolis\\
Parc Valrose, 06108 Nice\\France}
\email{cguillar@math.unice.fr}

\author{Leo Tzou}
\address{Department of Mathematics\\
Stanford University\\
Stanford, CA 94305, USA.}
\email{ltzou@math.stanford.edu.}

\begin{abstract}
On a fixed smooth compact Riemann surface with boundary $(M_0,g)$, we show that for the Schr\"odinger operator $\Delta +V$ with potential $V\in C^{1,\alpha}(M_0)$ for some $\alpha>0$, 
the Dirichlet-to-Neumann map $\mc{N}|_{\Gamma}$ measured on an open set $\Gamma\subset \pl M_0$ determines 
uniquely the potential $V$. We also discuss briefly the corresponding consequences for potential scattering at $0$ frequency
on Riemann surfaces with asymptotically Euclidean or asymptotically hyperbolic ends.
\end{abstract}

\maketitle

\begin{section}{Introduction}
The problem of determining the potential in the Schr\"odinger operator by boundary measurement  
goes back to Calder\'on \cite{Ca}. Mathematically, it amounts to ask if one can detect some data 
from boundary measurement in a domain (or manifold) $\Omega$ with boundary. The typical model to have in mind is the Schr\"odinger operator
$P:=\Delta_g+V$ where $g$ is a metric and $V$ a potential, then we define the Cauchy data space by 
\[\mc{C}:=\{(u|_{\pl\Omega},\pl_\nu u|_{\pl\Omega})\,; \,u \in H^1(\Omega), \ u \in \ker P\}\]
where $\pl_\nu$ is the interior pointing normal vector field to $\pl\Omega$.\\ 

The first natural question is the following \emph{full data} inverse problem: does the Cauchy data space determine uniquely the metric
$g$ and/or the potential $V$? 
In a sense, the most satisfying known results are when the domain $\Omega\subset \rr^n$ is already known and $g$ is
the Euclidean metric, then the recovery of $V$ has been proved in dimension $n>2$ by Sylvester-Uhlmann \cite{SU},
and very recently in dimension $2$ by Bukgheim \cite{Bu} when the domain is simply connected. A
related question is the conductivity problem 
which consists in taking $V=0$ and replacing $\Delta_g$ by $-{\rm div}\sigma\nabla$ where $\sigma$
is a definite positive symmetric tensor. An elementary observation shows that 
the problem of recovering an sufficiently smooth isotropic conductivity (i.e. $\sigma=\sigma_0{\rm Id}$ for a function $\sigma_0$) 
is contained in the problem above of recovering a potential $V$. 
For domain of $\rr^2$, Nachman \cite{N} used the $\bar{\pl}$ techniques to show 
that the Cauchy data space determines the conductivity. Recently a new approach developed by Astala and P\"aiv\"arinta in \cite{AP} improved this result to assuming that the conductivity is only a $L^\infty$ scalar function. This was later generalized to $L^\infty$ anisotropic conductivities by Astala-Lassas-P\"aiv\"arinta in \cite{ALP}. We notice that there still 
are rather few results in the direction of recovering the Riemannian manifold $(\Omega,g)$ 
when $V=0$, for instance the surface case by Lassas-Uhlmann \cite{LU} (see also \cite{Be,HM}), the real-analytic manifold case 
by Lassas-Taylor-Uhlmann \cite{LTU} (see also \cite{GSB} for the Einstein case),  the case of manifolds admitting limiting Carleman weights 
and in a same conformal class by Dos Santos Ferreira-Kenig-Salo-Uhlmann \cite{DKSU}.\\

The second natural, but harder, problem is the \emph{partial data} inverse problem: if $\Gamma_1$ and $\Gamma_2$ are open subsets of $\pl\Omega$, 
does the partial Cauchy data space for $P$
\[\mc{C}_{\Gamma_{1},\Gamma_2}:=\{(u|_{\Gamma_1},\pl_\nu u|_{\Gamma_2})\,; \,u\in H^1(\pl M_0), \, Pu=0,\,  u=0\,\, {\rm in }\,\,
\pl \Omega\setminus \Gamma_1\}\]
determine the domain $\Omega$, the metric, the potential? 
For a fixed domain of $\rr^n$, the recovery of the potential if $n>2$ with partial data measurements was initiated by Bukhgeim-Uhlmann \cite{BU} and later improved by 
Kenig-Sj\"ostrand-Uhlmann \cite{KSU} to the case where $\Gamma_1$ and $\Gamma_2$ are respectively open subsets of the "front" and "back" ends of the domain. We refer the reader to the references for a more precise formulation of the problem. In dimension $2$, the recent works of Imanuvilov-Uhlmann-Yamamoto \cite{IUY} solves the
problem for fixed domains $\Omega$ of $\rr^2$ in the case when $\Gamma_1 = \Gamma_2$ and when the potential are in $C^{2,\alpha}(\Omega)$ for some $\alpha>0$.\\ 

In this work, we address the same question when the background domain is a fixed Riemann surface with boundary.  
We prove the following recovery result: 
\begin{theorem}
\label{identif}
Let $(M_0,g)$ be a smooth compact Riemann surface with boundary and let $\Delta_g$ be its positive Laplacian. 
For $\alpha\in(0,1)$, let $V_1,V_2\in C^{1,\alpha}(M_0)$ be two real potentials and, for $i=1,2$, let 
\begin{equation}\label{cauchydata}
\mc{C}^\Gamma_i=:\{(u|_{\Gamma},\pl_\nu u|_{\Gamma})\,; \,u \in H^1(M_0), (\Delta_g+V_i)u=0,\,\,  u=0\,\,{\rm on }\,\, \pl M_0\setminus \Gamma\}
\end{equation}
be the respective Cauchy partial data spaces. If $\mc{C}^\Gamma_1=\mc{C}^\Gamma_2$ then $V_1=V_2$.
\end{theorem}
\end{section}
Here the space $C^{1,\alpha}(M_0)$ is the usual H\"older space for $\alpha\in(0,1)$.
Notice that when $\Delta_g+V_i$ do not have $L^2$ eigenvalues for the Dirichlet condition,
the statement above can be given in terms of Dirichlet-to-Neumann operators. 
Since $\Delta_{\hat{g}}=e^{-2\varphi}\Delta_g$ when $\hat{g}=e^{2\varphi}g$ for some function $\varphi$, it is clear that 
in the statement in Theorem \ref{identif}, we only need to fix the conformal class of $g$ instead of the metric $g$ 
(or equivalently to fix the complex structure on $M_0$). In particular, the smoothness assumption of the Riemann 
surface with boundary is not really essential since we can change it conformally to 
make it smooth and for the Cauchy data space, this just has the effect of changing the potential conformally (we only need this new potential to be $C^{1,\alpha}$).  
Observe also that Theorem \ref{identif} implies that, for a fixed Riemann surface with boundary $(M_0,g)$, the Dirichlet-to-Neumann map on $\Gamma$ for the operator $u\to -{\rm div}_{g}(\gamma \nabla^gu)$ determines the isotropic conductivity $\gamma$ if $\gamma\in C^{3,\alpha}(M_0)$ in the sense that two conductivities giving rise to the same Dirichlet-to-Neumann are equal. 
This is a standard observation by transforming the conductivity problem to a potential problem 
with potential $V:=(\Delta_g \gamma^{\demi})/\gamma^\demi$. So our result also extends that of Henkin-Michel \cite{HM} in the case
of isotropic conductivities.\\ 

The method to reconstruct the potential follows \cite{Bu,IUY} and is based on the construction of a large set of special complex geometric optic solutions of 
$(\Delta_g+V)u=0$, more precisely if $\Gamma_0= \pl M_0\setminus \Gamma$ is the set where we do not know the 
Dirichlet-to-Neumann operator, we construct solutions of the form 
$u={\rm Re}(e^{\Phi/h}(a+r(h))+e^{{\rm Re}(\varphi)/h}s(h)$ with $u|_{\Gamma_0}=0$ where $h>0$ is a small parameter, $\Phi$ and $a$ are holomorphic functions 
on $(M_0,g)$, independent of $h$, $||r(h)||_{L^2}=O(h)$ while $||s(h)||_{L^2}=O(h^{3/2}|\log h|)$  as $h\to 0$. The idea of \cite{Bu} to reconstruct $V(p)$ for $p\in M_0$
is to take $\Phi$ with a non-degenerate critical point at $p$ and then use stationary phase as $h\to 0$. In our setting, the function $\Phi$ needs to be purely real on
$\Gamma_0$ an Morse with a prescribed critical point at $p$. One of our main contribution is 
a geometric construction of the holomorphic Carleman weights $\Phi$ satisfying such conditions.
We should point out that we use a quite different method than in \cite{IUY} to construct this weight, and we believe that our method 
simplifies their construction even in their case. A Carleman estimate on the surface for this degenerate weight needs to be proved, and we follows ideas of \cite{IUY}. 
We manage to improve the regularity of the potential to $C^{1,\alpha}$ instead of $C^{2,\alpha}$ in \cite{IUY}. 
We finally notice that we provide a proof, in an Appendix, of the fact that the partial Cauchy data space $\mc{C}^\Gamma$
determine a potential $V\in C^{0,\alpha}(M_0)$ on $\Gamma$ (for $\alpha>0$).\\ 

In Section \ref{invscat}, we obtain two inverse scattering results as corollary of Theorem \ref{identif}:
first for partial data scattering at $0$ frequency for $\Delta+V$ on asymptotically hyperbolic surfaces with potential decaying at the conformal infinity, and secondly for full data scattering at $0$ frequenncy for $\Delta+V$ with $V$ compactly supported on an asymptotically Euclidean surface.\\

Another straightforward corollary in the asymptotically Euclidean case full data setting
is the recovery of a compactly supported potential from the scattering operator at a positive frequency. The proof
is essentially the same as for the operator $\Delta_{\rr^n}+V$ once we know Theorem \eqref{identif}, so we omit it. 

\begin{section}{Harmonic and Holomorphic Morse Functions on a Riemann Surface}

\subsection{Riemann surfaces}
We start by recalling few elementary definitions and results about Riemann surfaces, see for instance \cite{FK} for more details.
Let $(M_0,g_0)$ be a compact connected smooth Riemannian surface with boundary $\pl M_0$. The surface $M_0$ can be considered as 
a subset of a compact Riemannian surface $(M,g)$, for instance by taking the double of $M_0$ and extending smoothly 
the metric $g_0$ to $M$. 
The conformal class of $g$ on the closed surface $M$ induces a structure of closed Riemann surface, i.e. a closed surface equipped with a complex structure via holomorphic charts $z_\alpha:U_{\alpha}\to \cc$. The Hodge star operator $\star$ acts on the cotangent bundle $T^*M$, its eigenvalues are 
$\pm i$ and the respective eigenspace $T_{1,0}^*M:=\ker (\star+i{\rm Id})$ and $T_{0,1}^*M:=\ker(\star -i{\rm Id})$
are sub-bundle of the complexified cotangent bundle $\cc T^*M$ and the splitting $\cc T^*M=T^*_{1,0}M\oplus T_{0,1}^*M$ holds as complex vector spaces.
Since $\star$ is conformally invariant on $1$-forms on $M$, the complex structure depends only on the conformal class of $g$.
In holomorphic coordinates $z=x+iy$ in a chart $U_\alpha$,
one has $\star(udx+vdy)=-vdx+udy$ and
\[T_{1,0}^*M|_{U_\alpha}\simeq \cc dz ,\quad T_{0,1}^*M|_{U_\alpha}\simeq \cc d\bar{z}  \]
where $dz=dx+idy$ and $d\bar{z}=dx-idy$. We define the natural projections induced by the splitting of $\cc T^*M$ 
\[\pi_{1,0}:\cc T^*M\to T_{1,0}^*M ,\quad \pi_{0,1}: \cc T^*M\to T_{0,1}^*M.\]
The exterior derivative $d$ defines the De Rham complex $0\to \Lambda^0\to\Lambda^1\to \Lambda^2\to 0$ where $\Lambda^k:=\Lambda^kT^*M$
denotes the real bundle of $k$-forms on $M$. Let us denote $\cc\Lambda^k$ the complexification of $\Lambda^k$, then
the $\pl$ and $\bar{\pl}$ operators can be defined as differential operators 
$\pl: \cc\Lambda^0\to T^*_{1,0}M$ and $\bar{\pl}:\cc\Lambda_0\to T_{0,1}^*M$ by 
\begin{equation}\label{defddbar}
\pl f:= \pi_{1,0}df ,\quad \bar{\pl}:=\pi_{0,1}df,
\end{equation}
they satisfy $d=\pl+\bar{\pl}$ and are expressed in holomorphic coordinates by
\[\pl f=\pl_zf\, dz ,\quad \bar{\pl}f=\pl_{\bar{z}}f \, d\bar{z}.\]  
with $\pl_z:=\demi(\pl_x-i\pl_y)$ and $\pl_{\bar{z}}:=\demi(\pl_x+i\pl_y)$.
Similarly, one can define the $\pl$ and $\bar{\pl}$ operators from $\cc \Lambda^1$ to $\cc \Lambda^2$ by setting 
\[\pl (\omega_{1,0}+\omega_{0,1}):= d\omega_{0,1}, \quad \bar{\pl}(\omega_{1,0}+\omega_{0,1}):=d\omega_{1,0}\]
if $\omega_{0,1}\in T_{0,1}^*M$ and $\omega_{1,0}\in T_{1,0}^*M$.
In coordinates this is simply
\[\pl(udz+vd\bar{z})=\pl v\wedge d\bar{z},\quad \bar{\pl}(udz+vd\bar{z})=\bar{\pl}u\wedge d{z}.\]
There is a natural operator, the Laplacian acting on functions and defined by 
\[\Delta f:= -2i\star \bar{\pl}\pl f =d^*d \]
where $d^*$ is the adjoint of $d$ through the metric $g$ and $\star$ is the Hodge star operator mapping 
$\Lambda^2$ to $\Lambda^0$ and induced by $g$ as well.

\subsection{Maslov Index and Boundary value problem for the $\overline\partial$ Operator}
In this subsection we consider the setting where $M$ is an oriented Riemann surface with boundary $\partial M_0$ and $\Gamma\subset \partial M_0$ is an open subset and we let $\Gamma_0=\pl M_0\setminus \Gamma$ be its complement in $\pl M_0$. 
Following \cite{mcduff}, we adopt the following notations:
let $E\to M_0$ be a complex line bundle with complex structure $J : E\to E$ and let $D : C^{\infty}(M_0,E) \to C^\infty(M_0,T^*_{0,1}\otimes E)$ be a Cauchy-Riemann operator with smooth coefficients on $M_0$, acting on sections of the bundle $E$. Observe that in the case when $E=M_0\x \cc$ is the trivial line bundle with the natural complex structure on $M_0$, then $D$ can be taken as the operator $\overline\partial$ introduced in \eqref{defddbar}. 
For $q>1$, we define
\[D_F : W^{\ell,q}_F(M_0,E) \to W^{\ell-1,q}(M_0, T_{0,1}^*M_0\otimes E)\]
where $F\subset E\mid_{\pl M_0}$ is a totally real subbundle (i.e. a subbundle such that $JF \cap F$ is the zero section) 
and $D_F$ is the restriction of $D$ to the $L^q$-based Sobolev space will $\ell$ derivatives and boundary condition $F$
\[W^{\ell,q}_F(M_0,E) := \{\xi\in W^{\ell,q}(M_0,E) \mid \xi(\partial M_0) \subset F\}.\]
In this setting, we have the following boundary value Riemann-Roch theorem stated in \cite{mcduff}:
\begin{theorem}
\label{boundaryrr}
Let $E \to M_0$ be a complex line bundle over an oriented compact Riemann surface with boundary and $F\subset E\mid_{\pl M_0}$ be a totally real subbundle.  Let $D$ be a smooth Cauchy-Riemann operator on $E$ acting on $W^{\ell,q}(M_0,E)$ for some $q>1$ 	and $\ell\in\nn$. 
Then\\
1) The following operators are Fredholm 
\[D_F : W^{\ell,q}_F(M_0, E) \to W^{\ell-1,q}(M_0,T_{0,1}^*M_0\otimes E)\]
\[D_F^* : W^{\ell,q}_F(M_0, T_{0,1}^*M_0\otimes E) \to W^{\ell-1,q}(M_0, E).\]
2) The real Fredholm index of $D_F$ is given by
\[{\rm Ind}(D_F) = \chi(M_0) + \mu(E,F)\]
where $\chi(M_0)$ is the Euler characteristic of $M_0$ and $\mu(E,F)$ is the boundary Maslov index of the subbundle $F$.\\
3) If $\mu(E,F) <0$, then  $D_F$ is injective, while if  $\mu(E,F) + 2\chi(M_0) >0$ the operator $D_F$ is surjective.
\end{theorem}

In the case of a trivial bundle $E=M_0\x \cc$ and $\pl M_0=\sqcup_{j=1}^m \pl_i M_0$ is a disjoint union of $m$
circles, the Maslov index can be defined (see \cite[p.554-555]{mcduff})  to be 
the degree of the map $\rho \circ \Lambda:\pl M_0\to \pl M_0$ where 
\[\Lambda|_{\pl_i M_0}: S^1\simeq \pl_i M_0\to{\rm GL}(1,\cc)/{\rm GL}(1,\rr)\] 
is the natural map given by the totally real subbundle (the space ${\rm GL}(1,\cc)/{\rm GL}(1,\rr)$ being the space of totally real subbundles of $\cc$ over a circle) and 
\[\rho: {\rm GL}(1,\cc)/{\rm GL}(1,\rr) \to S^1, \quad \rho(A.{\rm GL}(1,\rr)):=\det(A^2)/\det(A^*A).\] 

As an application, we obtain the following (here and in what follows, $H^m(M_0):=W^{m,2}(M_0)$):
\begin{coro}\label{surject}
(i) For $q>1$ and $k\in \nn_0$, let $\omega \in W^{k,q}(M_0,T^{*}_{0,1}M_0)$, then there exists $u\in W^{k+1,q}(M_0)$ holomorphic on $M_0$, 
real valued on $\Gamma_0$, such that $\bar{\pl}u=\omega$.\\   
(ii) For $m> 1/2$, let $f\in H^{m}(\pl M_0)$ be a real valued function, then there exists a holomorphic function $v\in H^{m+\demi}(M_0)$ such that
${\rm Re}(v)|_{\Gamma_0}=f$.\\
(iii) For $k\in\nn$ and $q>1$, the space of $W^{k,q}(M_0)$ holomorphic functions on $M_0$ which are real valued on $\Gamma_0$ is infinite dimensional.
\end{coro}
\noindent{\bf Proof}. (i) Let $L\in\nn$ be arbitrary large and let us identify the boundary as a disjoint union of circles $\partial M_0= \coprod_{i = 1}^m \pl_iM_0$ where each $\pl_iM_0 \simeq S^1$. Since $\Gamma$ will be the piece of the boundary where we know the Cauchy data space, it is sufficient to assume that $\Gamma$ is a connected non-empty open segment of $\pl_1M_0 = S^1$, and which can thus be defined in a coordinate $\theta$ (respecting the orientation of the boundary) by 
$\Gamma = \{\theta \in S^1 \mid 0< \theta< 2\pi/k\}$ for some integer $k$.
Define the totally real subbundle of $F\subset E|_{\partial M_0} = \coprod_{j = 1}^l (\pl_jM_0 \times \C)$ by the following: on 
$\pl_1M_0 \simeq S^1$ parametrized by $\theta \in [0,2\pi]$, define $F_\theta = e^{ia(\theta)}\R$. Where $a: [0,2\pi] \to \R$ is a smooth nondecreasing function such that $a(\theta) = 0$ in a neighbourhood $[0,\eps]$ of $0$, $a(2\pi/k) = 2L\pi$ 
for some $L\in\nn$, and $a(\theta) = 2L\pi$ for all $\theta > 2\pi/k$. For the rest of $\pl_2M_0,.., \pl_m M_0$, we  just let 
$F|_{\pl_iM_0} = S^1 \times \R$. The Maslov index $\mu(E,F)$ is then given by $2L$
and so, by theorem \ref{boundaryrr}, $D_F$ is surjective if $2\chi(M_0)+2L>0$. Since $L$ can be taken 
as large as we want this achieves the proof of (i).\\ 

(ii) Let $w\in H^{m+\demi}(M_0)$ be a real function with boundary value $f$ on $\pl M_0$, then by (i) there exists $R\in H^{m+1/2}(M_0)$ 
such that $i\bar{\pl}R=-\bar{\pl}w$ and $R$ purely real on $\Gamma_0$, thus $v:=iR+w$ is holomorphic such that ${\rm Re}(v)=f$ on $\Gamma_0$.\\ 

(iii) Taking the subbundle $F$ as in the proof of (i), we have that $\dim \ker D_F=\chi(M_0)+2L$ if $L$ satisfies $2\chi(M_0)+2L>0$, 
and since $L$ can be taken as large as we like-, this concludes the proof. 
\qed

\begin{lemma}
\label{control the zero}
Let $\{p_0, p_1,..,p_n\}\subset M_0$ be a set of $n+1$ disjoint points. Let $c_1,\dots, c_K\in\cc$, $N\in\nn$, 
and let $z$ be a complex coordinate near $p_0$ such that $p_0=\{z=0\}$. 
Then if $p_0\in{\rm int}(M_0)$, there exists a holomorphic function $f$ on $M_0$ with zeros of order at least $N$ at each $p_j$,
such that $f$ is real on $\Gamma_0$ and $f(z)=c_0 + c_1z +...+ c_K z^K+ O(|z|^{K+1})$ in the coordinate $z$. If 
$p_0\in \pl M_0$, the same is true except that $f$ is not necessarily real on $\Gamma_0$.
\end{lemma}
\noindent{\bf Proof}.
First, using linear combinations and induction on $K$, it suffices to prove the Lemma for any $K$ and $c_0=\dots=c_{K-1}=0$, which we now show.  
Consider the subbundle $F$ as in the proof of (i) in Corollary \ref{surject}. 
The Maslov index $\mu(E,F)$ is  given by $2L$ and so for each $N\in\nn$, one can take $L$ large enough to have $\mu(F,E) + 2\chi(M_0) \geq 2N(1+n)$. Therefore by Theorem \ref{boundaryrr} the dimension of the kernel of $\overline\partial_F$ will be greater than $2(n+1)N$. 
Now, since for each $p_j$ and complex coordinate $z_j$ near $p_j$, 
the map $u\to (u(p_j),\pl_{z_j}u(p_j),\dots,\pl_{z_j}^{N-1}u(p_j))\in \cc^{N}$ is linear, 
this implies that there exists a non-zero element $u\in \ker D_F$ which has zeros of order at least $N$ 
at all $p_j$.

First, assume that $p_0\in {\rm int}(M_0)$ and we want the desired Taylor expansion at $p_0$ in the coordinate $z$.  
In the coordinate $z$, one has $u(z) = \alpha z^M+O(|z|^{M+1})$ for some $\alpha\not=0$ and $M\geq N$. 
Define the function $r_K(z) = \chi(z)\frac{c_K}{\alpha} z^{-M+K}$ 
where $\chi(z)$ is a smooth cut-off function supported near $p_0$ and
which is $1$ near $p_0=\{z=0\}$. Since $M \geq N> 1$, this function has a pole at $p_0$ and trivially extends smoothly to 
$M_0\backslash \{p_0\}$, which we still call $r_K$. Observe that the function is holomorphic in a neighbourhood of $p_0$ but not at $p_0$ where it is only meromorphic, so that in $M_0\setminus \{p_0\}$, $\overline\partial r_K$ is a smooth and compactly supported section of $T^*_{0,1}M_0$ 
and therefore trivially extends smoothly to $M_0$ (by setting its value to be $0$ at $p_0$) to a one form denoted $\omega_K$. 
By the surjectivity assertion in Corollary \ref{surject}, 
there exists a smooth function $R_K$ satisfying $\overline\partial R_K = - \omega_K$ and that $R_K|_{\Gamma_0} \in \R$. 
We now have that $R_K + r_K$ is a holomorphic function on $M\backslash\{p_0\}$ 
meromorphic with a pole of order $M-K$ at $p_0$, and in coordinate $z$ one has $z^{M-K}(R_K(z)+r_K(z))= c_K+O(|z|)$. 
Setting $f_K= u(R_K + r_K)$, we have the desired holomorphic function. Note that $f$ also vanish to order $N$ at all 
$p_1,\dots,p_n$ since $u$ does. This achieves the proof.

Now, if $p_0\in \pl M_0$, we can consider a slightly larger smooth domain of $M$ containing $M_0$ and we apply the the result above. 
\qed

\subsection{Morse holomorphic functions with prescribed critical points} 
The main result of this section is the following 
\begin{proposition}
\label{criticalpoints}
Let $p$ be an interior point of $M_0$ and $\epsilon >0$ small. 
Then there exists a holomorphic function $\Phi$ on $M_0$ which is Morse on $M_0$ (up to the boundary) and real valued on $\Gamma_0$, 
which has a critical point $p'$ at distance less than $\eps$ from $p$ and such that 
${\rm Im}(\Phi(p'))\not=0$. 
\end{proposition}
Let $\mc{O}$ be a connected open set of $M$ such that $\bar{\mc{O}}$ is a smooth surface with boundary, with $M_0\subset \bar{\mc{O}}$ and $\Gamma_0\subset\pl\bar{\mc{O}}$.
Fix $k>2$ a large integer, we denote by $C^k(\bar{\mc{O}})$ the Banach space of $C^k$ real valued functions on $\bar{\mc{O}}$. 
Then the set of harmonic functions on $\bar{\mc{O}}$ which are in the Banach space
$C^{k}(\bar{\mc{O}})$ (and smooth in $\mc{O}$ by elliptic regularity) 
is the kernel of the continuous map $\Delta:C^k(\bar{\mc{O}})\to C^{k-2}(\bar{\mc{O}})$, 
and so it is a Banach subspace of $C^k(\bar{\mc{O}})$. 
The set $\mc{H}\subset C^k(\bar{\mc{O}})$ of harmonic functions $u$ in $C^k(\bar{\mc{O}})$ such there exists 
$v\in C^{k}(\bar{\mc{O}})$ harmonic with $u+iv$ holomorphic on $\mc{O}$ is a Banach subspace of $C^{k}(\bar{\mc{O}})$ of 
finite codimension. Indeed, let $\{\gamma_1,..,\gamma_N\}$ be a homology basis for $\mc{O}$, then
\[\mc{H}=\ker L , \textrm{ with } L: \ker\Delta\cap C^k(\bar{\mc{O}})\to \cc^N \textrm{ defined by }
L(u):=\Big(\frac{1}{\pi i}\int_{\gamma_j}\pl u\Big)_{j=1,\dots,N}.\]
For all $\Gamma_0'\subset \pl M_0$ such that the complement of $\Gamma_0'$ contains an open subset, 
we define 
\[\mc{H}_{\Gamma_0'} := \{u \in\mc{H} ; u|_{\Gamma_0'} =0\}.\]
We now show
\begin{lemma}\label{morsedense}
The set of functions $u\in\mc{H}_{\Gamma_0'}$ which are Morse in $\mc{O}$ 
is residual (i.e. a countable intersection of open dense sets) 
in $\mc{H}_{\Gamma_0'}$ with respect to the $C^k(\bar{\mc{O}})$ topology.
\end{lemma}
\noindent{\bf Proof}. We use an argument very similar to those used by Uhlenbeck \cite{Uh}.
We start by defining $m: \mc{O}\times \mc{H}_{\Gamma_0'}\to T^*\mc{O}$ by $(p,u) \mapsto (p,du(p))\in T_p^*\mc{O}$. 
This is clearly a smooth map, linear in the second variable, moreover $m_u:=m(.,u)=(\cdot, du(\cdot))$ is  
Fredholm since $\mc{O}$ is finite dimensional. The map $u$ is a Morse function if and only if 
$m_u$ is transverse to the zero section, denoted $T_0^*\mc{O}$, of $T^*\mc{O}$, ie. if 
\[\textrm{Image}(D_{p}m_u)+T_{m_u(p)}(T_0^*\mc{O})=T_{m_u(p)}(T^*\mc{O}),\quad \forall p\in \mc{O} \textrm{ such that }m_u(p)=(p,0),\]
which is equivalent to the fact that the Hessian of $u$ at critical points is 
non-degenerate (see for instance Lemma 2.8 of \cite{Uh}). 
We recall the following transversality theorem (\cite[Th.2]{Uh}):
\begin{theorem}\label{transversality}
Let $m : X\times \mc{H}_{\Gamma_0'} \to W$ be a $C^k$ map, where $X$, $\mc{H}_{\Gamma_0'}$, and $W$ are separable Banach manifolds 
with $W$ and $X$ of finite dimension. Let $W'\subset W$ be a submanifold such that $k>\max(1,\dim X-\dim W+\dim W')$.
If $m$ is transverse to $W'$  then the set 
$\{u\in \mc{H}_{\Gamma_0'}; m_u \textrm{ is transverse to } W'\}$ is dense in $\mc{H}_{\Gamma_0'}$, more precisely 
it is a residual set. 
\end{theorem}
We want to apply it with $X:=\mc{O}$, $W:=T^*\mc{O}$ and $W':=T^*_0\mc{O}$, and the map $m$ is defined above. 
We have thus proved Lemma \ref{morsedense} if one can show that $m$ is transverse to $W'$. 
Let $(p,u)$ such that $m(p,u)=(p,0)\in W'$. Then identifying $T_{(p,0)}(T^*\mc{O})$ with $T_p\mc{O}\oplus T^*_p\mc{O}$, one has
\[D_{(p,u)}m(z,v)=(z,dv(p)+{\rm Hess}_p(u)z)\]
where ${\rm Hess}_pu$ is the Hessian of $u$ at the point $p$, viewed as a linear map from $T_p\mc{O}$ to $T^*_p\mc{O}$. 
To prove that $m$ is transverse to $W'$
we need to show that $(z,v)\to (z, dv(p)+{\rm Hess}_p(u)z)$ is onto from $T_p\mc{O}\oplus \mc{H}_{\Gamma_0'}$ to $T_p\mc{O}\oplus T^*_p\mc{O}$, which 
is realized for instance if the map $v\to dv(p)$ from $\mc{H}_{\Gamma_0'}$ to  $T_p^*\mc{O}$ is onto.
But from Lemma \ref{control the zero}, we know that there exist holomorphic functions $v$ and $\tilde v$ on $\mc{O}$  such that $v$ and $\tilde v$ are purely real on $\Gamma_0'$. Clearly the imaginary parts of $v$ and $\tilde v$ belong to $\mc{H}_{\Gamma_0'}$. Furthermore, for a given complex coordinate $z$ near $p=\{z=0\}$, we can arrange them to have series expansion $v(z)=z+ O(|z|^2)$ and ${\tilde v}(z)=iz + O(|z|^2)$ around the point $p$. We see, by coordinate computation of the exterior derivative of ${\rm Im}(v)$ and ${\rm Im}(\til{v})$, that 
 $d\, {\rm Im}(v)(p)$ and $d \,{\rm Im}(\tilde v)(p)$ are linearly independent at the point $p$. This shows 
 our claim and ends the proof of Lemma \ref{morsedense} by using Theorem \ref{transversality}.\qed\\
We now proceed to show that the set of all functions $u\in \mc{H}_{\Gamma_0'}$ such that $u$ has no degenerate critical points on $\Gamma_0'$ is also residual.

\begin{lemma}
\label{nonvanishing normal derivative}
For all $p\in \Gamma_0'$ and $k\in\nn$, there exists a holomorphic function $u\in C^{k}(\bar{\mc{O}})$, such that ${\rm Im}(u)|_{\Gamma_0'}=0$  and 
$\pl u(p) \neq 0$.
\end{lemma}
\noindent{\bf Proof}. The proof is quite similar to that of Lemma \ref{control the zero}. By Lemma \ref{control the zero}, we can choose a holomorphic function $v\in C^k(\bar{\mc{O}})$ such that $v(p)=0$ and ${\rm Im}(v)|_{\Gamma_0'}=0$, then either $\pl v(p)\not=0$ and we are done, or $\pl v(p)=0$. 
Assume now the second case and let $M\in\nn$ be the order of $p$ as a zero of $v$. By Riemann mapping theorem, there is a conformal
mapping from a neighbourhood $U_p$ of $p$ in $\bar{\mc{O}}$ to a neighbourhood $\{|z|<\eps, {\rm Im}(z)\geq 0\}$ of the 
real line ${\rm Im}(z)=0$ in $\cc$, and one can assume that $p=\{z=0\}$ in these complex coordinates. 
Take $r(z)=\chi(z)z^{-M+1}$ where $\chi\in C_0^\infty(|z|\leq \eps)$ is a real valued function with $\chi(z)=1$ in $\{|z|<\eps/2\}$. Then 
$\bar{\pl} r$ 
vanishes in the pointed disc $0<|z|<\eps/2$ and it is a compactly supported smooth section of $T^*_{1,0}\bar{\mc{O}}$ outside, it can thus be 
extended trivially to a smooth section of $T^*_{1,0}\bar{\mc{O}}$ denoted by $\omega$. 
We can then use (i) of Corollary \ref{surject}:
there is a function $R$ such that $\bar{\pl} R=-\omega$ and ${\rm Im}(R)|_{\Gamma_0'}=0$, and
so $\bar{\pl}(R+r)=0$ in $\mc{O}\setminus \{p\}$ and $R+r$ is real valued on $\Gamma_0'$ (remark that $r$ is real valued on $\Gamma_0'$)
and has a pole at $p$ of order exactly $M-1$. We conclude that $u:=v(R+r)$ satisfies the desired properties, it 
vanishes at $p$ but with non zero complex derivative at $p$. 
\qed

\begin{lemma}
\label{degeneracy condition}
Let $\Gamma'_0\subset \partial \mc{O}$ be an open set of the boundary. Let $\phi : \mc{O} \to \R$ be a harmonic function with $\phi|_{\Gamma_0'} = 0$. Let $p\in\Gamma_0'$ be a critical point of $\phi$, then it is nondegenerate 
if and only if $\partial_\tau\partial_\nu u \neq 0$ where $\pl_\tau$ and $\pl_\nu$ denote respectively the tangential and normal 
derivatives along the boundary.
\end{lemma}
\noindent{\bf Proof}. By Riemman mapping theorem, there is a conformal transformation mapping a neighbourhood of $p$ in $\bar{\mc{O}}$
to a half-disc $D:=\{|z|<\eps,{\rm Im}(z)\geq 0\}$ and $\pl\bar{\mc{O}}=\{{\rm Im}(z)=0\}$ near $p$. Denoting $z=x+iy$, one has
$(\pl_x^2+\pl_y^2)\phi=0$ in $D$ and $\pl^2_x\phi|_{y=0}=0$, which implies $\pl_y^2\phi(p)=0$. Since $\pl_\nu=e^{f}\pl_y$ and 
$\pl_\tau=e^{f}\pl_x$ for some smooth function $f$, and since $d\phi(p)=0$, the conclusion is then straightforward.
\qed\\
Let $N^{*}\pl \bar{\mc{O}}$ be the conormal-bundle of $\pl\bar{\mc{O}}$ and $N^*\Gamma_0'$ be the restriction of this bundle to $\Gamma_0'$. Denote the zero sections of these bundles respectively by $N^*_0\pl\bar{ \mc{O}}$ and $N^{*}_0\Gamma_0'$. 
We now define the map
\[b : \Gamma_0' \times \mc{H}_{\Gamma_0'}  \to N^{*}\Gamma_0', \quad b(p,u):=(p,\pl_\nu u).\] 
For a fixed $u\in \mc{H}_{\Gamma_0'}$, we also define  $b_u(\cdot) := b(\cdot, u)$. Simple computations yield the
\begin{lemma}
\label{mixed partial non-vanishing}
Suppose that $p\in\Gamma_0'$ is such that $\partial_\nu u(p) = 0$, then $\partial_\tau\partial_\nu u(p)\neq 0$ if and only if
\[{\rm Image}(D_pb_u) + T_{(p,0)}(N_0^{*}\Gamma_0') = T_{(p,0)}(N^{*}\Gamma_0').\]
\end{lemma}
\noindent{\bf Proof}. This can be seen by the fact that for all $p\in \Gamma'_0$ such that $b_u(p) = (p,0)$, 
\[D_pb_u : T_p\Gamma_0'\to T_{(p,0)}(N^{*}\Gamma_0') \simeq T_p\Gamma_0' \oplus N^{*}_p\Gamma_0'\] 
is given by $w\mapsto (w, \partial_\tau\partial_\nu u(p)w)$.\qed\\

At a point $(p,u)$ such that $b(p,u) = 0$, a simple computation yields that the differential $D_{(p,u)}b : T_p\Gamma_0' \times \mc{H}_{\Gamma_0'} \to T_{(p,\pl_\nu u(p))}(N^{*}\Gamma_0')$ is given by $(w,u') \mapsto (w, \partial_\tau\partial_\nu u(p)w + \partial_\nu u'(p))$. This observation combined with Lemma \ref{nonvanishing normal derivative} shows that for all $(p,u) \in \Gamma_0'\times \mc{H}_{\Gamma_0'}$ such that $b(p,u) = (p,0)$, $b$ is transverse to $N^{*}_0\Gamma_0'$ at $(p,0)$. Now we can apply Theorem \ref{transversality} with $X = \Gamma_0'$, $W = N^{*}\Gamma_0'$ and $W' = N^{*}_0\Gamma_0'$ we see that the set $\{u\in \mc{H}_{\Gamma_0'}; b_u \textrm{ is transverse to }N^{*}_0\Gamma_0'\}$ is residual in $\mc{H}_{\Gamma_0'}$. In view of Lemmas \ref{degeneracy condition}, we deduce the 
\begin{lemma}
\label{boundarymorse}
The set of functions $u\in \mc{H}_{\Gamma_0'}$ such that $u$ has no degenerate critical point on $\Gamma_0'$ is residual in $\mc{H}_{\Gamma_0'}$.
\end{lemma}
Observing the general fact that finite intersection of residual sets remains residual, the combination of Lemma \ref{boundarymorse} and Lemma \ref{morsedense} yields
\begin{coro}
\label{boundarymorse dense}
The set of functions $u\in\mc{H}_{\Gamma_0'}$ which are Morse in $\mc{O}$ and have no degenerate critical points on $\Gamma_0'$ 
is residual in $\mc{H}_{\Gamma_0'}$ with respect to the $C^k(\bar{\mc{O}})$ topology. In particular, it is dense.
\end{coro}
We are now in a position to give a proof of the main proposition of this section.\\\\
\noindent{\bf Proof of Proposition \ref{criticalpoints}}. As explained above, choose $\mc{O}$ in such a way that $\bar{\mc{O}}$ is a smooth surface with boundary, containing $M_0$, that $\Gamma_0\subset\partial \mc{O}$ and $\mc{O}$ contains $\partial M_0 \backslash \overline{\Gamma_0}$. Let $\Gamma_0'$ be an open subset of the boundary of $\bar{\mc{O}}$ such that the closure of $\Gamma_0$ is contained in $\Gamma_0'$ and $\partial\bar{\mc{O}}\backslash\overline{\Gamma_0'}\not=\emptyset$. Let $p$ be an interior point of $M_0$. By lemma \ref{control the zero}, there exists a holomorphic function $f = u + i v$ on $\bar{\mc{O}}$ such that $f$ is purely real on $\Gamma_0'$, $v(p)= 1$, and $df(p) = 0$ (thus $v\in \mc{H}_{\Gamma_0'}$).\\

By Corollary \ref{boundarymorse dense}, there exist a sequence $(v_j)_j$ of Morse functions $v_j \in \mc{H}_{\Gamma_0'} $ 
such that $v_j\to v$ in $C^k(M_0)$ for any fixed $k$ large. By Cauchy integral formula, there exist harmonic conjugates $u_j$ of $v_j$ such that $f_j := u_j + i v_j \to f$ in $C^k(M_0)$. Let $\eps>0$ be small and let $U\subset\mc{O}$ 
be a neighbourhood containing $p$ and no other critical points of $f$,
and with boundary a smooth circle of radius $\eps$. In complex local coordinates near $p$, 
we can identify $\pl f$ and $\pl f_j$ to holomorphic functions on an open set of $\cc$.
Then by Rouche's theorem, it is clear that $\pl{f_j}$ has precisely one zero in $U$ and $v_j$ never vanishes in $U$ if $j$ is large enough.\\\\
Fix $\Phi$ to be one of the $f_j$ for $j$ large enough. By construction, $\Phi$ is Morse in $\mc{O}$ and has no degenerate critical points on $\overline{\Gamma_0} \subset\Gamma_0'$. 
We notice that, since the imaginary part of $\Phi$ vanishes on all of $\Gamma_0'$, it is clear from
the reflection principle applied after using the Riemann mapping theorem (as in the proof of Lemma \ref{degeneracy condition})
that no point on $\overline{\Gamma_0}\subset \Gamma_0'$ can be an accumulation point for critical points. Now $\partial M_0 \backslash \overline{\Gamma_0}$ is contained in the interior of $\mc{O}$ and therefore no points on  $\partial M_0 \backslash \overline{\Gamma_0}$ can be an accumulation point of critical points. Since $\Phi$ is Morse in the interior of $\mc{O}$, there are no degenerate critical points on $\partial M_0 \backslash \overline{\Gamma_0}$. This ends the proof.\qed

\end{section}

\begin{section}{Carleman Estimate for Harmonic Weights with Critical Points}

In this section, we prove a Carleman estimate using harmonic weight with non-degenerate critical points, in 
way similar to \cite{IUY}. Let us define $\Gamma_0:=\{p\in \pl M_0; \pl_\nu \varphi(p)=0\}$ and let $\Gamma:=\pl M_0\setminus \Gamma_0$ its 
complement.

\begin{proposition}
\label{carlemanestimate}
Let $(M_0,g)$ be a smooth Riemann surface with boundary, and let $\varphi: M_0\to\R$ be a $C^k(M_0)$ harmonic Morse function for $k$ large.  Then for all $V\in L^\infty(M_0)$ there exists an $h_0>0$ such that for all $h\in(0,h_0)$ and $u\in C^\infty(M)$ with $u|_{\pl M_0}=0$, we have 
\begin{equation}
\begin{gathered}
\label{carleman}
\frac{1}{h}\|u\|^2_{L^2(M_0)} + \frac{1}{h^2}\|u |d\varphi|\|^2_{L^2(M_0)} + \|du\|^2_{L^2(M_0)}+
\|\pl_\nu u\|^2_{L^2(\Gamma_0)} \\
\leq C\Big(\|e^{-\varphi/h}(\Delta_g + V)e^{\varphi/h} u\|_{L^2(M_0)}^2+\frac{1}{h}\|\pl_\nu u\|^2_{L^2(\Gamma)}\Big)
\end{gathered}
\end{equation}
where $\pl_\nu$ is the exterior unit normal vector field to $\pl M_0$. 
\end{proposition}
\noindent\textbf{Proof}.  
We start by modifying the weight as follows:
if $\varphi_0:=\varphi: M_0 \to \R$ is a real valued harmonic Morse function with critical points $\{p_1,\dots,p_N\}$ in the interior of $M_0$, we 
let $\varphi_j:M_0 \to \R$ be harmonic functions such that $p_j$ is not a critical point of $\varphi_j$ for $j = 1,\dots,N$, their existence 
is insured by  Lemma \ref{control the zero}.
For all $\epsilon >0$, we define the convexified weight $\varphi_{\eps} := \varphi - \frac{h}{2\epsilon}(\sum_{j = 0}^N|\varphi_j|^2)$.\\

To prove the estimate, we shall localize in charts $\Omega_j$ covering the surfaces. These charts will be taken so that if 
$\Omega_j\cap \pl M_0\not=\emptyset$, then $\Omega_j\cap \pl M_0\simeq S^1$ is a connected component of $\pl M_0$. Moreover, by Riemann
mapping theorem (e.g. Lemma 3.2 of \cite{MT}), this chart can be taken to be a neighbourhood of $|z|=1$ in $\{z\in \cc; |z|\leq 1\}$ and such that the metric 
$g$ is conformal to the Euclidean metric $|dz|^2$.
\begin{lemma}
Let $\Omega$ be a chart of $M_0$ as above and $\varphi_{\eps}:\Omega \to \R$ be as above. Then there are constants $C,C'>0$ such that for all $\omega\in C^\infty(M)$ supported in $\Omega$ and $h>0$ small enough, the following estimate holds:
\begin{equation}
\label{dbar carleman}
\frac{C}{\epsilon}\|\omega\|_{L^2(M_0)}^2 +C'\Big(-{\rm Im}(\cjg\pl_\tau \omega,\omega\cjd_{L^2(\pl M_0)})
+\frac{1}{h}\int_{\pl M_0}|\omega|^2\pl_\nu\varphi_\eps {\rm dv}_g\Big)\leq \|e^{-\varphi_{\eps}/h}\bar{\pl} e^{\varphi_{\eps}/h}
\omega\|^2_{L^2(M_0)}
\end{equation} 
where $\pl_{\nu}$ and $\pl_\tau$ denote respectively the exterior pointing normal vector fields and 
its rotation by an angle $+\pi/2$.
\end{lemma}
\noindent{\bf Proof}.
We use complex coordinates $z=x+iy$ in the chart $\Omega$ where $u$ is supported. Observe that the Lebesgue measure 
$dxdy$ is bounded below and above by ${\rm dv}_g$, $g$ is conformal to $|dz|^2$ and the boundary terms
in \eqref{dbar carleman} depend only on the conformal class, 
it suffices to prove the estimates with respect to $dxdy$ and the Euclidean metric. 
We thus integrate by parts with respect to $dxdy$ 
\begin{equation}\label{carl1}
\begin{split}
4\|e^{-\varphi_{\eps}/h}\bar{\pl}e^{\varphi_{\eps}/h}u\|^2 =& 
\|(\partial_{x} + \frac{i\pl_y\varphi_\eps}{h})u + (i\partial_{y} + \frac{\pl_x\varphi_\eps}{h})u\|^2\\
=& \|(\partial_{x} + \frac{i\pl_y\varphi_\eps}{h})u\|^2 + \|(i\partial_{y} + \frac{\pl_x\varphi_\eps}{h})u\|^2\\ 
&+ \frac{2}{h} \int_{\Omega}\Big(\Delta\varphi_\eps|u|^2 -\demi\pl_x\varphi_\eps.\pl_x|u|^2-\demi\pl_y\varphi_\eps.\pl_y|u|^2\Big)\\
&+ \frac{2}{h}\int_{\pl M_0}\pl_{\nu}\varphi_\eps |u|^2
-2\int_{M_0} \Big(\pl_x{\rm Re}(u).\pl_y{\rm Im}(u)-\pl_x{\rm Im}(u).\pl_y{\rm Re}(u)\Big)\\
=&\|(\partial_{x} + \frac{i\pl_y\varphi_\eps}{h})u\|^2 + \|(i\partial_{y} + \frac{\pl_x\varphi_\eps}{h})u\|^2
+\frac{1}{h} \int_{\Omega}\Delta\varphi_\eps|u|^2\\ 
&+\frac{1}{h}\int_{\pl M_0}\pl_{\nu}\varphi_\eps |u|^2+2\int_{\pl M_0}\pl_\tau {\rm Re}(u).{\rm Im}(u).
\end{split}
\end{equation}
where $\Delta:=-(\pl_x^2+\pl_y^2)$, $\pl_\nu$ is the exterior pointing normal vector field to the boundary 
and $\pl_\tau$ is the tangent vector field to the boundary (i.e. $\pl_\nu$ rotated with an angle $\pi/2$) for the Euclidean metric
$|dz|^2$. 
Then $\cjg u\Delta \varphi_\eps,u\cjd = \frac{h}{\epsilon}(|d\varphi_0|^2 + |d\varphi_1|^2 + .. + |d\varphi_N|^2)|u|^2$,
since $\varphi_j$ are harmonic, so the proof follows from the fact that $|d\varphi_0|^2+|d\varphi_1|^2 + .. + |d\varphi_N|^2$ is uniformly 
bounded away from zero.\qed\\

The main step to go from (\ref{dbar carleman}) to (\ref{carleman}) is the following lemma which is a slight modification of the proof in \cite{IUY}:
\begin{lemma}\label{ce on patch}
With the same assumptions as in Proposition \ref{carlemanestimate}, and if $\Omega$ is either an interior chart of $(M,g)$ 
or a chart containing a whole boundary connected component, then there are positive constants $c$ and $C$ such that for all $\eps>0$ small, all $0<h\ll \eps$ and all smooth function $u$ supported in $\Omega$ with $u|_{\pl M_0}=0$, we have
\begin{equation}\label{estimateDelta}
\begin{gathered}
C\Big(\|e^{-\varphi_{\eps}/h}\Delta_g e^{\varphi_{\eps}/h} u\|_{L^2(M_0)}^2 + 
\frac{1}{h}\|\pl_\nu u\|^2_{L^2(\Gamma)}\Big)\geq \\
\frac{c}{\epsilon}\Big(\frac{1}{h}\|u\|_{L^2(M)}^2 + \frac{1}{h^2}\|u |d\varphi|\|_{L^2(M)}^2 +  \frac{1}{h^2}\|u |d\varphi_{\eps}|\|_{L^2(M)}^2+ \|du\|_{L^2(M)}^2\Big) +\|\pl_\nu u\|_{L^2(\pl M_0)}^2 
\end{gathered}
\end{equation}
where $\pl_\nu$  denote the unit normal vector field to $\pl M_0$.
\end{lemma}
\noindent{\bf Proof}. Since the norms induced by the metric $g$ in the chart are conformal 
to Euclidean norms, and there is $f$ smooth such that $\Delta_g=-e^{2f}(\pl_x^2+\pl_y^2)=e^{2f}\Delta$ in the complex coordinate chart, it suffices to get the estimate \eqref{estimateDelta} for Euclidean norms and Laplacian.
Clearly, we can assume $u\in H^1_0(M)$ to be real valued without loss of generality.  Now let $Q(z)$ be a holomorphic
function in $\Omega$ which has no zeros, then by (\ref{carl1}) we have 
\[\begin{gathered}
\|Qe^{-\varphi_{\eps}/h}\Delta e^{\varphi_{\eps}/h} u\|^2  = 16\|e^{-\varphi_{\eps}/h}\bar\partial e^{\varphi_{\eps}/h}Qe^{-\varphi_{\eps}/h}\partial e^{\varphi_{\eps}/h} u\|^2 \geq \\
\frac{C}{\epsilon}\|\partial u + \frac{\pl\varphi_{\eps}}{h} u\|^2  -4{\rm Im}\Big(\int_{\pl M_0}\pl_\tau \omega.\bar{\omega}\Big)
+\frac{4}{h}\int_{\pl M_0}|\omega|^2\pl_\nu\varphi_\eps  
\end{gathered}\]
with $\omega:=Qe^{-\varphi_{\eps}/h}\partial e^{\varphi_{\eps}/h} u$ and here everything is measured 
with respect to Euclidean metric $dx^2+dy^2$ and measure $dxdy$. 
Since $u|_{\pl M_0}=0$, one has  
$\omega|_{\pl M_0}=(A+iB)\pl_\nu u$ where $A+iB=Q(\cjg\pl_\nu,\pl_x\cjd-i\cjg\pl_\nu,\pl_y\cjd)$ and 
$-{\rm Im}(\pl_\tau \omega.\bar{\omega})=(\pl_\tau A.B-A\pl_\tau B)|\pl_\nu u|^2$
we deduce that for some $c>0$ 
\begin{equation}
\begin{gathered}\label{boundaryterms}
\|Qe^{-\varphi_\eps/h}\Delta e^{\varphi_\eps/h} u\|^2 \geq\\
 \frac{c}{\eps}\Big(\|du\|^2 + \frac{1}{h^2} \|u |d\varphi_{\eps}|\|^2 + \frac{2}{h}\cjg\pl_xu, u\pl_x\varphi_\eps\cjd+ 
\frac{2}{h} \cjg\pl_y u ,u\pl_y \varphi_{\eps}\cjd\Big)\\
+4\int_{\pl M_0}(A\pl_\tau B-\pl_\tau AB)|\pl_\nu u|^2
+\frac{4}{h}\int_{\pl M_0}|Q|^2|\pl_\nu u|^2\pl_\nu\varphi_\eps. 
\end{gathered}
\end{equation}
Using the fact that $u$ is real valued, that $\varphi$ is harmonic and that $\sum_{j=0}^N|d\varphi_j|^2$ is uniformly 
bounded away from $0$, we see that
\[
\frac{2}{h}\cjg\pl_xu, u\pl_x\varphi_\eps\cjd+ \frac{2}{h} \cjg\pl_y u,u\pl_y \varphi_\eps\cjd = 
\frac{1}{h}\cjg u,u \Delta\varphi_\eps\cjd \geq \frac{C}{\epsilon}\|u\|^2
\]
for some $C>0$ and therefore,
\begin{equation}\label{ineq1b}
\|Qe^{-\varphi_{\eps}/h}\Delta e^{\varphi_{\eps}/h} u\|^2 \geq \frac{c}{\epsilon}(\|du\|^2 + \frac{1}{h^2} \|u |d\varphi_{\eps}|\|^2  + \frac{C}{\epsilon}\|u\|^2)+ \textrm{ boundary terms}.
\end{equation}
Now if the diameter of the support of $u$ is chosen small (with size depending only on $|{\rm Hess}\varphi_0|(p)$) 
with a unique critical point $p$ of $\varphi_0$ inside, 
one can use integration by parts and the fact that the critical point is non-degenerate to obtain 
\begin{equation}\label{ineq2} 
\|\bar{\pl} u\|^2 + \frac{1}{h^2} \|u |\pl \varphi_0|\|^2 \geq \frac{1}{h}\left|\int \pl_{\bar{z}}(u^2)\bbar{\pl_z\varphi_0} dxdy \right|
\geq \frac{1}{h}\left|\int u^2\,\bbar{\pl_z^2 \varphi_0}\, dxdy \right|\geq
\frac{C'}{h}\|u\|^2
\end{equation}
for some $C'>0$. Clearly the same estimate holds trivially if $\Omega$ does not contain critical point of $\varphi_0$.  
Using a partition of unity $(\theta_j)_j$ in $\Omega$ and absorbing terms of the form $||u\bar{\pl}\theta_j||^2$ into the right hand side,
one obtains \eqref{ineq2} for any function $u$ supported in $\Omega$ and vanishing at the boundary. 
Thus, combining with \eqref{ineq1b}, there are positive constants $c,c',C''$ such that for $h$ small enough
\[\begin{gathered}
\frac{c}{\epsilon}(\|du\|^2 + \frac{1}{h^2} \|u |d\varphi_{\eps}|\|^2  + \frac{C}{\epsilon}\|u\|^2) 
\geq \frac{c}{\epsilon}(\|du\|^2 + \frac{1}{h^2} \|u |d\varphi_0|\|^2-\frac{C''}{\eps^2}\|u\|^2) 
\\
\geq\frac{c'}{\epsilon}(\|du\|^2 + \frac{1}{h^2} \|u |d\varphi_0|\|^2 + \frac{1}{h}\|u\|^2 ). 
\end{gathered}\]
Combining now with \eqref{ineq1b}  and using that $|Q|$ is bounded below gives
\[\|e^{-\varphi_{\eps}/h}\Delta e^{\varphi_{\eps}/h} u\|^2 \geq \frac{c'}{\epsilon}(\|du\|^2 + \frac{1}{h^2} \|u |d\varphi|\|^2  + \frac{1}{h}\|u\|^2)+ \textrm{ boundary terms}.\]
Let us now discuss the boundary terms in \eqref{boundaryterms}. If $\varphi_j$ are taken so that $\pl_\nu\varphi_j=0$ on $\Gamma_0$, then 
$\pl_\nu\varphi_\eps=0$ on $\Gamma_0$ and $\pl_\nu\varphi_\eps=\pl_\nu\varphi+O(h/\eps)$ on $\Gamma$ and thus
\[\frac{1}{h}\int_{\pl M_0}|Q|^2|\pl_\nu u|^2|\pl_\nu\varphi_\eps|\leq \frac{C_1}{h}\int_{\Gamma}|\pl_\nu u|^2\]
for some constant $C_1$. We finally claim that there exist $Q$ with no zeros in $\Omega$ such that $A\pl_\tau B-B\pl_\tau A$ is bounded below 
by a positive constant on $\pl M_0\cap \Omega$. Indeed, since the chart near a connected component can be taken to be an interior 
neighbourhood of the circle $|z|=1$ in $\cc$, one can take $A+iB=e^{it}$ where $t\in S^1$ parametrize the boundary
component, so that $A\pl_\tau B-B\pl_\tau A=1$ since $\pl_\tau=\pl_t$ for the Euclidean metric. Since moreover
$\cjg\pl_\nu,\pl_x\cjd-i\cjg\pl_\nu,\pl_y\cjd=A-iB=e^{-it}$, we deduce that on the boundary $Q(t)=e^{2it}$ and so it suffices
to take $Q(z)=z^2$.  This achieves the proof.
\qed\\

\noindent{\bf Proof of Proposition \ref{carlemanestimate}}. Using triangular inequality and absorbing the term $||Vu||^2$
into the left hand side of \eqref{carleman}, it suffices to prove \eqref{carleman} with  $\Delta_g$ instead of $\Delta_g+V$.
Let $v \in C_0^\infty(M)$, we have by Lemma \ref{ce on patch} that there exist constants $c,c',C,C'>0$ such that
\[\begin{gathered}
\frac{c}{\epsilon}\Big(\frac{1}{h}\|e^{-\varphi_{\eps}/h}v\|^2 + \frac{1}{h^2}\|e^{-\varphi_{\eps}/h}v |d\varphi|\|^2 
+\frac{1}{h^2}\|e^{-\varphi_{\eps}/h}v |d\varphi_\eps|\|^2+ \|e^{-\varphi_{\eps}/h}dv\|^2\Big) +\|e^{-\varphi_\eps/h}\pl_\nu v\|^2_{\Gamma_0} \\
\leq \sum_j \frac{c'}{\epsilon}\Big(\frac{1}{h}\|e^{-\varphi_{\eps}/h}\chi_jv\|^2 + \frac{1}{h^2}\|e^{-\varphi_{\eps}/h}\chi_jv |d\varphi|\|^2 
  +\frac{1}{h^2}\|e^{-\varphi_{\eps}/h}\chi_jv |d\varphi_\eps|\|^2\\
  + \|e^{-\varphi_{\eps}/h}d(\chi_jv)\|^2\Big)
  + \| e^{-\varphi_\eps/h}\pl_\nu v\|_{\Gamma_0} \\ 
\leq  C\Big(\sum_j \|e^{-\varphi_{\eps}/h}\Delta_g  (\chi_jv)\|^2+\|e^{-\varphi_\eps/h}\pl_\nu v\|^2_{\Gamma}\Big)\\ \leq C'\Big(\|e^{-\varphi_{\eps}/h}\Delta_g  v\|^2 +  \|e^{-\varphi_{\eps}/h} v\|^2 + \|e^{-\varphi_{\eps}/h} dv\|^2 +\|e^{-\varphi_\eps/h}\pl_\nu v\|^2_{\Gamma}\Big)
\end{gathered}
\]
where $(\chi_j)_j$ is a partition of unity associated to the complex charts $\Omega_j$ on $M$. 
Since constants on both sides are independent of $\eps$ and $h$, we can take $\eps$ small enough so that 
$C' \|e^{-\varphi_\eps/h} v\|^2 + C' \|e^{-\varphi_\eps/h} dv\|^2$ can be absorbed to the left side. 
Now set $v = e^{\varphi_\eps/h} w$ with $w|_{\pl M_0}=0$, then we have 
\[\begin{gathered}
\frac{1}{h}\|w\|^2 + \frac{1}{h^2}\|w |d\varphi|\|^2 +\frac{1}{h^2}\|w |d\varphi_{\eps}|\|^2 +\|dw\|^2 +\|\pl_\nu \omega\|^2_{\Gamma_0}\\
\leq C\Big(\|e^{-\varphi_{\eps}/h}\Delta_g e^{\varphi_{\eps}/h} w\|^2 +\|\pl_\nu \omega\|_{\Gamma}^2\Big)
\end{gathered}\]
Finally, fix $\epsilon >0$ and set $u: = e^{\frac{1}{\epsilon}\sum_{j = 0}^N|\varphi_j|^2}w$ and use the fact that $e^{\frac{1}{\epsilon}\sum_{j = 0}^N|\varphi_j|^2}$ is independent of $h$ and bounded uniformly away from zero 
and above, we then obtain the desired estimate for $0<h \ll\epsilon$.\qed

\end{section}

\begin{section}{Complex Geometric Optics on a Riemann Surface}\label{CGOriemann}
As explained in the Introduction, the method for identifying the potential at a point $p$ is to construct complex geometric optic solutions depending on a small parameter $h>0$, with phase a Carleman weight (here a Morse holomorphic function), and such that the phase has a non-degenerate critical point at $p$, in order to apply the stationary phase method. In this section, the potential $V$ has the regularity $V\in C^{1,\alpha}(M_0)$ for some $\alpha>0$.
 
Choose $p\in {\rm int}(M_0)$ such that there exists a holomorphic function $\Phi=\varphi+i\psi$ which is Morse on $M_0$, $C^k$ in $M_0$ for large $k\in\nn$ and such that $\partial\Phi(p) = 0$ and $\Phi$ has only finitely many critical points in $M_0$. Furthermore we ask that $\Phi$ is purely real on $\Gamma_0$. By Proposition \ref{criticalpoints} such points $p$ form a dense subset of $M_0$. Given such a holomorphic function, the purpose of this section is to construct solutions $u$ on $M_0$ of $(\Delta +V)u = 0$ of the form
\begin{equation}
\label{cgo}
u = e^{\Phi/h}(a + ha_0 + r_1) + \bbar{e^{\Phi/h}(a + ha_0 +r_1)} + e^{\varphi/h}r_2\ \textrm{ with } \ \ u|_{\Gamma_0} = 0
\end{equation}
for $h>0$ small, where $a$ is holomorphic and $u\in C^k(M_0)$ for large $k\in\nn$, $a_0\in H^{2}(M_0)$ is holomorphic, moreover $a(p)\not=0$ and $a$ vanishes to high order at all other critical points $ p' \in M_0$ of $\Phi$. Furthermore, we ask that the holomorphic function $a$ is purely imaginary on $\Gamma_0$. The existence of such a holomorphic function is a consequence of Lemma \ref{control the zero}. Given such a holomorphic function on $M_0$ we consider a compactly supported extension to $M$, still denoted $a$.\\

The remainder terms $r_1,r_2$ will be controlled as $h\to 0$ and have particular properties near the critical points of $\Phi$.
More precisely, $r_2$ will be a $O_{L^2}(h^{3/2}|\log h|)$ and $r_1$ will be of the form $h\til{r}_{12}+o_{L^2}(h)$ where $\til{r}_{12}$ is independent of $h$, 
which can be used to obtain sufficient informations from the stationary phase method in the identification process.

\begin{subsection}{Construction of $r_1$}\label{constr1}
We shall  construct $r_1$ to satisfy
\[e^{-\Phi/h}(\Delta_g +V)e^{\Phi/h}(a + r_1) = O_{L^2}(h|\log h|)\]
and $r_1 = r_{11} + h r_{12}$.
We let $G$ be the Green operator of the Laplacian on the smooth surface with boundary 
$M_0$ with Dirichlet condition, so that $\Delta_gG={\rm Id}$ on $L^2(M_0)$. In particular this implies 
that $\bar{\partial}\partial G=\frac{i}{2}\star^{-1}$ where $\star^{-1}$ is the inverse of $\star$ mapping functions to $2$-forms.
We extend $a$ to be a compactly supported $C^k$ function on $M_0$ and we will search for $r_{1}\in H^{2}(M_0)$ satisfying
$||r_1||_{L^2}=O(h)$ and
\begin{equation}
\label{dequation}
e^{-2i\psi/h}\partial e^{2i\psi/h} r_1 = -\pl G (aV) + \omega + O_{H^1}(h|\log h|)
\end{equation}
where $\omega$ is a smooth holomorphic 1-form on $M_0$. 
Indeed, using the fact that $\Phi$ is holomorphic we have
\[e^{-\Phi/h}\Delta_ge^{\Phi /h}=-2i\star  \bar{\pl} e^{-\Phi/h}\pl e^{\Phi/h}=-2i\star  \bar{\pl} e^{-\frac{1}{h}(\Phi-\bar{\Phi})}\pl e^{\frac{1}{h}(\Phi-\bar{\Phi})}=
-2i\star \bar{\pl}e^{-2i\psi/h}\pl e^{2i\psi/h}\]
and applying $-2i\star\bar{\pl}$ to \eqref{dequation}, we obtain (note that $\pl G(aV)\in C^{2,\alpha}(M_0)$ by elliptic regularity)
\[e^{-\Phi/h}(\Delta_g+V)e^{\Phi/h}r_1=-aV+O_{L^2}(h|\log h|).\]
We will choose $\omega$ to be a smooth holomorphic $1$-form on $M_0$ such that at all 
critical point $p'$ of $\Phi$ in $M_0$, the form $b:=\pl G(aV) - \omega$  with value in $T^*_{1,0}M_0$ vanish to the highest possible order.
Writing $b=b(z)dz$  in local complex coordinates, $b(z)$ is $C^{2+\alpha}$ by elliptic regularity and 
we have $-2i\pl_{\bar{z}}b(z)=aV$, therefore $\pl_z\pl_{\bar{z}}b(p')=\pl^2_{\bar{z}}b(p')=0$ at each critical point $p'\not=p$ by construction 
of the function $a$.  Therefore, we deduce that at each critical point $p'\neq p$, $\pl G(aV)$ has Taylor series expansion $\sum_{j = 0}^2 c_j z^j + O(|z|^{2+\alpha})$ for some $N$ large. That is, all the lower order terms of the Taylor expansion of $\pl G(aV)$ around $p'$ are polynomials of $z$ only. 
\begin{lemma}\label{formomega}
Let $\{p_0,...,p_N\}$ be finitely many points on $M_0$ and let $\theta$ be a $C^{2,\alpha}$ section of $T^*_{1,0}M_0$.
Then there exists a $C^k$ holomorphic function $f$ on $M$ with $k\in\nn$ large, such that $\omega = \partial f$ satisfies the following: in complex local coordinates $z$ near 
$p_j$ , one has $\pl_z^\ell\theta(p_j)=\pl_z^\ell\omega(p_j)$ for $\ell=0,1,2$, where $\theta=\theta(z)dz$ and $\omega=\omega(z)dz$.
\end{lemma}
\noindent{\bf Proof}. This is a direct consequence of Lemma \ref{control the zero}.
\qed\\
Applying this to the form $\pl G(aV)$ and using the observation we made above, we can construct a $C^k$ holomorphic form $\omega$ such that in local coordinates $z$ centered at a critical point 
$p'$ of $\Phi$ (i.e $p'=\{z=0\}$ in this coordinate), we have for $b=\pl G(aV)-\omega=b(z)dz$
\begin{equation}\label{decayofb}
\begin{gathered}
|\pl_{\bar{z}}^m\pl^{\ell}_z b(z)|=O(|z|^{2+\alpha-\ell-m}) , \textrm{ for } \ell+m\leq 2 ,\quad \textrm{ if }p'\not=p\\
|b(z)|=O(|z|) , \quad \textrm{ if }p'=p.  
 \end{gathered}
 \end{equation}

Now, we let $\chi_1\in C_0^\infty(M_0)$ be a cutoff function supported in a small neighbourhood $U_p$ of the critical point $p$ and identically $1$ near $p$, and 
$\chi\in C_0^\infty(M_0)$ is defined similarly with $\chi =1$ on the support of $\chi_1$.
We will construct $r_1= r_{11} +h r_{12}$ in two steps : first, we will construct $r_{11}$ to solve equation \eqref{dequation} 
locally near the critical point $p$ of $\Phi$ and then we will 
construct the global correction term $r_{12}$ away from $p$ by using the extra vanishing of $b$ in \eqref{decayofb} at the other critical points.

We define locally in complex coordinates centered at $p$ and containing the support of $\chi$
\[r_{11}:=\chi e^{-2i\psi/h}R(e^{2i\psi/h}\chi_1b)\] 
where $Rf(z) := -(2\pi i)^{-1}\int_{\R^2} \frac{1}{\bar{z}-\bar{\xi}}f d\bar{\xi}\wedge d\xi$ for $f\in L^\infty$ compactly supported 
is the classical Cauchy-Riemann operator inverting locally $\pl_z$ ($r_{11}$ is extended by $0$ outside the neighbourhood of $p$).
The function $r_{11}$ is in $C^{3+\alpha}(M_0)$ and we have 
\begin{equation}\begin{gathered}\label{r11}
e^{-2i\psi/h}\pl(e^{2i\psi/h}r_{11}) = \chi_1(-\pl G(aV) + \omega) + \eta\\
\textrm{ with }\eta:= e^{-2i\psi/h}R(e^{2i\psi/h}\chi_1b)\pl\chi.
\end{gathered}
\end{equation}
We then construct $r_{12}$ by observing that  $b$ vanishes to order $2+\alpha$ at critical points of $\Phi$ other than $p$ (from \eqref{decayofb}), and 
$\pl \chi=0$ in a neighbourhood of any critical point of $\psi$, so we can find $r_{12}$ satisfying
\[2ir_{12}\pl\psi = (1-\chi_1)b .\] 
This is possible since both $\pl\psi$ and the right hand side are valued in $T^*_{1,0}M_0$, $\pl \psi$ has finitely many isolated $0$ on $M_0$: 
$r_{12}$ is then a function which is in $C^{2,\alpha}(M_0\setminus{P})$ where $P:=\{p_1,\dots, p_N\}$ is the set of critical points other than $p$,
it extends to a $C^{1,\alpha}(M_0)$  and it satisfies in local complex coordinates $z$ near each $p_j$ 
\[ |\pl_{\bar{z}}^\beta\pl_z^\gamma r_{12}(z)|\leq C|z-p_j|^{1+\alpha-\beta-\gamma} , \quad \beta+\gamma\leq 2.\]
by using also the fact that $\pl \psi$ can be locally be considered as holomorphic function with a zero of order $1$ at each $p_j$.
This implies that $r_1\in H^2(M_0)$ and  we have 
\[ e^{-2i\psi/h}\pl(e^{2i\psi/h}r_1) = b+h\pl r_{12}=-\pl G(aV)-\omega+ h\pl r_{12} + \eta.\]
Now the first error term $||\pl r_{12}||_{H^1(M_0)}$ is bounded by 
\[||\pl r_{12}||_{H^1(M_0)}\leq C\left( \left|\left|\frac{(1-\chi_1)b(z)}{\pl_z \psi(z)}\right|\right|_{H^2(U_p)}\right)\leq
C \]
for some constant $C$, where we used the fact that $\frac{(1-\chi_1)b(z)}{\pl_z \psi(z)}$ is in $H^2(U_p)$ and independent of $h$.
To deal with the $\eta$ term, we need the following 
\begin{lemma}\label{termeta}
The following estimates hold true 
\[\begin{gathered} 
||\eta||_{H^2}=O(|\log h|),\quad \|\eta\|_{H^1}\leq O(h|\log h|),  \quad ||r_{1}||_{L^2}=O(h), \quad ||r_1-h\til{r}_{12}||_{L^2}=o(h) 
\end{gathered} \]
where $\til{r}_{12}$ solves $2i\til{r}_{12}\pl\psi = b$.
\end{lemma}
\textbf{Proof}. We start by observing that 
\begin{equation}\label{decomposition} 
\begin{gathered}
||r_1||_{L^2}\leq \left| \left| \chi e^{-2i\psi/h}R (e^{2i\psi/h}\chi_1 b)-h\frac{\chi_1b}{2\pl_z\psi}\right|\right|_{L^2(U_p)}+
h||\til{r}_{12}||_{L^2(M_0)}+ h\left|\left|\frac{\eta(z)}{2\pl_z\psi}\right|\right|_{L^2(U_p)},\\
||r_1-h\til{r}_{12}||_{L^2}\leq \left| \left| \chi e^{-2i\psi/h}R (e^{2i\psi/h}\chi_1 b)-h\frac{\chi_1b}{2\pl_z\psi}\right|\right|_{L^2(U_p)}+
 h\left|\left|\frac{\eta(z)}{2\pl_z\psi}\right|\right|_{L^2(U_p)}
\end{gathered}
\end{equation}
and we will show that the first and last terms in the right hand sides are $o(h)$ while  $h||\til{r}_{12}||_{L^2(M_0)}$ is $O(h)$.
The first term is estimated in Proposition 2.7 of \cite{IUY}, it is a $o(h)$, while the last term is clearly bounded by
$Ch||\eta||_{L^2}$ and the middle one by $Ch$ for some constant $C$ by using that $\pl \psi$ does not vanish on the support of $\eta$ and the fact that
$b$ vanishes at critical points of $\psi$. 
Now are going to estimate the $H^2$ norms of $\eta$. Locally in complex coordinates $z$ centered at $p$ (ie. $p=\{z=0\}$), we have 
\begin{equation} \label{definitioneta}
\eta(z)= -\pl_z\chi(z)e^{-\frac{2i\psi(z)}{h}}\int_{\cc} e^{\frac{2i\psi(\xi)}{h}}\frac{1}{\bar{z}-\bar{\xi}} \chi_1(\xi)b(\xi) \frac{d\xi_1d\xi_2}{\pi} , 
\quad \xi=\xi_1+i\xi_2.
\end{equation} 
Since $b$ is $C^{2,\alpha}$ in $U$, we decompose $b(\xi)=\cjg\nabla b(0),\xi\cjd+\til{b}(\xi)$ using Taylor formula, so we have 
$\til{b}(0)=\pl_\xi \til{b}(0)=0$ and we split the integral \eqref{definitioneta} with $\cjg\nabla b(0),\xi\cjd$ and $\til{b}(\xi)$.  
Since the integrand with the $\cjg\nabla b(0),\xi\cjd$ is smooth and compactly supported in $\xi$ (recall that $\chi_1=0$ on the support of $\pl_z\chi$),
we can apply stationary phase to get that 
\[\left|\pl_z\chi(z)e^{-\frac{2i\psi(z)}{h}}\int_{\cc} e^{\frac{2i\psi(\xi)}{h}}\frac{1}{\bar{z}-\bar{\xi}} \chi_1(\xi)\cjg\nabla b(0),\xi\cjd \frac{d\xi_1d\xi_2}{\pi}\right|\leq Ch^2\]
uniformly in $z$.
Now set $\til{b}_z(\xi)=\pl_z\chi(z)\chi_1(\xi)\til{b}(\xi)/(\bbar{z-\xi})$ which is $C^{2,\alpha}$ in $\xi$ and smooth in $z$. 
Let $\theta\in C_0^\infty([0,1))$ be a cutoff function which is equal to $1$ near $0$ and 
set $\theta_h(\xi):=\theta(|\xi|/h)$, then we have by integrating by parts 
\begin{equation}\label{intbyparts}
\begin{split}
\int_{\cc} e^{\frac{2i\psi(\xi)}{h}}\til{b}_z(\xi) d\xi_1d\xi_2= &
h^2\int_{{\rm supp}(\chi_1)} e^{\frac{2i\psi(\xi)}{h}}\pl_{\bar{\xi}}\left(\frac{1-\theta_h(\xi)}{2i\pl_{\bar{\xi}} \psi}
\pl_{\xi}\left(\frac{\til{b}_z(\xi)}{2i\pl_\xi\psi}\right)\right) d\xi_1d\xi_2\\
&-h\int_{{\rm supp}(\chi_1)} e^{\frac{2i\psi(\xi)}{h}}\theta_h(\xi)\pl_{\xi}\left(\frac{\til{b}_z(\xi)}{2i\pl_\xi\psi}\right) d\xi_1d\xi_2 .
\end{split}\end{equation}
Using polar coordinates with the fact that $\til{b}_z(0)=0$, 
it is easy to check that the second term in \eqref{intbyparts} is bounded uniformy in $z$ 
by $Ch^{2}$. 
To deal with the first term, we use $\til{b}_z(0)=\pl_\xi \til{b}_z(0)=\pl_{\bar{\xi}}\til{b}_z(0)=0$ and 
a straightforward computation in polar coordinates shows
that the first term of \eqref{intbyparts}  is bounded uniformly in $z$ by $Ch^{2}|\log(h)|$.
We conclude that
\[ ||\eta||_{L^2}\leq C||\eta||_{L^\infty}\leq Ch^{2}|\log h|.\]
It is also direct to see that the same estimates holds with a loss of $h^{-2}$ for any derivatives in $z,\bar{z}$ of order less or equal to $2$, 
since they only hit the $\chi(z)$ factor, the $(\bar{z}-\bar{\xi})^{-1}$ factor or the oscillating term $e^{-2i\psi(z)/h}$.
So  we deduce that 
\[ ||\eta||_{H^2}= O(|\log h|).\]
and this ends the proof.
\qed\\

We summarize the result of this section with the following
\begin{lemma}
\label{solve equation to next order}
Let $k\in\nn$ be large and $\Phi\in C^k(M_0)$ be a holomorphic function on $M_0$ which is Morse in $M_0$ with a critical point at $p\in {\rm int}(M_0)$. 
Let $a\in C^k(M_0)$ be a holomorphic function on $M_0$ vanishing to high order at any critical point of $\Phi$ other than $p$. 
Then there exists $r_1 \in H^{2}(M_0)$ such that $||r_1||_{L^2}=O(h)$ and
\[e^{-\Phi/h}(\Delta +V)e^{\Phi/h}(a + r_1) =O_{L^2}(h|\log h|).\]
\end{lemma}
\end{subsection}
\begin{subsection}{Construction of $a_0$}
We have constructed the correction terms $r_1$ which solves the Schr\"odinger equation to order $h$ as stated in Lemma \ref{solve equation to next order}. In this subsection, we will construct  a holomorphic function $a_0$ which annihilates the boundary value of the solution on 
$\Gamma_0$. In particular, we have the following
\begin{lemma}
\label{a0 and a1}
There exists a holomorphic function $a_0\in H^2(M_0)$ independent of $h$ such that 
\[e^{-\Phi/h}(\Delta +V)e^{\Phi/h}(a + r_1 + h a_0) = O_{L^2}(h|\log h|)\]
and 
\[[e^{\Phi/h}(a + r_1 + h a_0)+\bbar{e^{\Phi/h}(a + r_1 + ha_0 )}]|_{\Gamma_0} =0.\]
\end{lemma}
\noindent{\bf Proof}. First, notice that $h^{-1}r_1|_{\pl M_0}=\til{r}_{12}|_{\pl M_0}\in H^{3/2}(\pl M_0)$ is independent of $h$.
Since $\Phi$ is purely real on $\Gamma_0$ and $a$ is purely imaginary on $\Gamma_0$, we see that this 
Lemma amounts to construct a holomorphic function $a_0\in H^{2}(M_0)$  with the boundary condition
\[{\rm Re}(\til{r}_{12}) + {\rm Re}(a_0)= 0 \, \textrm{ on }\,\Gamma_0.\]
To construct $a_0$, it suffices to use (ii) in Corollary \ref{surject}.  
\qed
\end{subsection}
\begin{subsection}{Construction of $r_2$}
The goal of this section is to complete the construction of the complex geometric optic solutions by the following proposition:
\begin{proposition}
\label{completecgo}
There exist solutions to $(\Delta +V)u = 0$ with boundary condition $u|_{\Gamma_0} = 0$ of the form \eqref{cgo} with $r_1$, $a_0$ constructed in the previous sections and $r_2$ satisfying $\|r_2\|_{L^2}= O(h^{3/2}|\log h|)$.
\end{proposition}
This is a consequence of the following Lemma (which follows from the Carleman estimate obtained above): 
\begin{lemma}
\label{standardargument}
If $V\in L^\infty(M_0)$ and $f\in L^2(M_0)$, then for all $h>0$ small enough, there exists a solution $v\in L^2$ to the boundary value problem
\[e^{\varphi/h}(\Delta_g +V) e^{-\varphi/h}v = f, \quad v|_{\Gamma_0} = 0,\]
satisfying the estimate
\[\|v\|_{L^2} \leq Ch^\demi\|f\|_{L^2}.\]
\end{lemma}
\noindent{\bf Proof}. The proof is the same as Proposition 2.2 of \cite{IUY}, we repeat the argument for the convenience of the reader.
Define for all $h >0$ the real vector space $\mc{A}:=\{u\in H_0^1(M_0); (\Delta_g+V)u\in L^2(M_0),\  \pl_\nu u\mid_{\Gamma}= 0\}$
equipped with the real scalar product 
\[(u,w)_{\mc{A}}:=\int_{M_0}e^{-2\varphi/h}(\Delta_gu+Vu)(\Delta_gw+Vw)\,{\rm dv}_g.\] 
Observe that since $\psi$ is constant along $\Gamma_0$, $\pl_\nu \varphi = 0$ on $\Gamma_0$. Therefore, we may apply the Carleman estimate of Proposition \eqref{carlemanestimate} to the weight $\varphi$ to assert that the space $\mc{A}$ is a Hilbert space equipped with the scalar product above. By using the same estimate, the linear functional $L:w\to \int_{M_0}e^{-\varphi/h}fw \,d\rm{v}_g$ on $\mc{A}$ is continuous and its  norm is bounded by $h^\demi||f||_{L^2}$.  By Riesz theorem, 
there is an element $u\in\mc{A}$ such that $(.,u)_{\mc{A}}=L$   and with norm bounded by the norm of 
$L$. It remains to take $v:=e^{-\varphi/h}(\Delta_g u+Vu)$ which solves $(\Delta_g+V)e^{-\varphi}v=e^{-\varphi/h}f$  and 
which in addition satisfies the desired norm estimate.
Furthermore, since
\[\int_{M_0} e^{-\varphi/h} v (\Delta_g + V) w {\rm dv}_g= \int_{M_0}e^{-\varphi/h}fw\, \rm{dv}_g\]
for all $w\in \mc{A}$, we have by Green's theorem
\[\int_{\pl M_0}e^{-\varphi/h}v\pl_\nu w {\rm dv}_g= 0=\int_{\Gamma_0}e^{-\varphi/h}v\pl_\nu w{\rm dv}_g \]
for all $w\in \mc{A}$. This implies $v = 0$ on $\Gamma_0$.
\qed\\

\noindent{\bf Proof of Proposition \ref{completecgo}}. We note that 
\[(\Delta +V)(e^{\Phi/h}(a + r_1 + ha_0) + \bbar{e^{\Phi/h}(a +  r_1 + ha_0)} + e^{\varphi/h}r_2 ) = 0\]
if and only if
\[e^{-\varphi/h}(\Delta + V) e^{\varphi/h} r_2=-e^{-\varphi/h}(\Delta +V)(e^{\Phi/h}(a + r_1 + ha_0) + \bbar{e^{\Phi/h}(a +  r_1 + h a_0)}).\]
By Lemma \ref{a0 and a1}, the right hand side of the above equation is $O_{L^2}(h|\log h|)$. Therefore, using Lemma \ref{standardargument} one can find such $r_2$ which satisfies
\[\|r_2\|_{L^2} \leq Ch^{3/2}|\log h|,\quad\quad r_2\mid_{\Gamma_0} = 0.\]
Since the ansatz $e^{\Phi/h}(a + r_1 + ha_0) + \bbar{e^{\Phi/h}(a + r_1 + ha_0)}$ is arranged to vanish on $\Gamma_0$, the solution 
\[u =e^{\Phi/h}(a + r_1 + ha_0) + \bbar{e^{\Phi/h}(a + r_1 + h a_0)} + e^{\varphi/h}r_2 \]
vanishes on $\Gamma_0$ as well.\qed
\end{subsection}

\end{section}
\begin{section}{Identifying the potential}
We now assume that $V_1,V_2\in C^{1,\alpha}(M_0)$ are two real valued potentials, with $\alpha>0$, such that the respective Cauchy data spaces
$\mc{C}^\Gamma_1,\mc{C}^\Gamma_2$ for the operators $\Delta_g+V_1$ and $\Delta_g+V_2$ on $\Gamma\subset \pl M_0$ are equal.
Let $\Gamma_0=\pl M_0\setminus \Gamma$ be the complement of $\Gamma$ in $\pl M_0$, and possibly by taking $\Gamma$ slightly smaller, we may assume that $\Gamma_0$ contains an open set.
Let $p\in M_0$ be an interior point of $M_0$ such that, using Proposition \ref{criticalpoints}, we can choose a holomorphic Morse function $\Phi=\varphi+i\psi$ on $M_0$ with $\Phi$ purely real on $\Gamma_0$, $C^k$ in $M_0$ for some large $k\in\nn$, with a critical point at $p$. Note that Proposition \ref{criticalpoints} states that we can choose $\Phi$ such that none of its critical points on the boundary are degenerate and such that critical points do not accumulate on the boundary.
\begin{proposition}
\label{identcritpts}
If the Cauchy data spaces agree, i.e. if $\mc{C}^{\Gamma}_1 = \mc{C}^{\Gamma}_2$ , then $V_1(p)= V_2(p)$.
\end{proposition}
\noindent{\bf Proof}. 
Let $a$ be a holomorphic function on $M_0$ which is purely imaginary on $\Gamma_0$ with $a(p)\neq 0$ and $a(p') = 0$ to large order for all other critical point $p'$ of $\Phi$. 
The existence of $a$ is insured by Lemma \ref{control the zero}.
Let $u_1$ and $u_2$ be $H^2$ solutions on $M_0$ to 
\[(\Delta_g +V_j)u_j = 0\]
constructed in Section \ref{CGOriemann} with $\Phi = \phi + i\psi$ for Carleman weight for $u_1$ and $-\Phi$ for $u_2$, thus of the form
\[u_1 = e^{\Phi/h}(a + h a_0 + r_1) + \bbar{e^{\Phi/h}(a + ha_0 + r_1)} + e^{\varphi/h} r_2 \]
\[u_2= e^{-\Phi/h}(a + h b_0 + s_1) + \bbar{e^{-\Phi/h}(a + h b_0 + s_1)} + e^{-\varphi/h} s_2 \]
and with boundary value $u_j|_{\pl M_0}=f_j$, where $f_j$ vanishes on $\Gamma_0$.
We can write by Green formula
\[\begin{split}
\int_{M_0}u_1(V_1 - V_2) {u_2} {\rm dv}_g&=- \int_{M_0} (\Delta_gu_1.  {u_2} - u_1.\Delta_g {u_2}){\rm dv}_g\\
&=-\int_{\pl M_0} (\pl_\nu u_1. {f_2}-f_1. \pl_\nu  {u_2}) {\rm dv}_{g}.
 \end{split}\]
Since the Cauchy data for $\Delta_g + V_1$ agrees on $\Gamma$ with that of $\Delta_g + V_2$, there exists a solution $v$ of the boundary value problem
\[(\Delta_g + V_2) v = 0,\ \ \ \ \ v|_{\pl M_0}= f_1,\]
satisfying $\pl_\nu v = \pl_\nu u_1$ on $\Gamma$. Since $f_j=0$ on $\Gamma_0$, this implies that 
\begin{equation}\label{integralid}\begin{split}
\int_{M_0}u_1(V_1 - V_2)  {u_2} {\rm dv}_g&=- \int_{M_0} (\Delta_gu_1.  {u_2} - u_1.\Delta_g {u_2}){\rm dv}_g=-\int_{\pl M_0} (\pl_\nu u_1 .{f_2}-f_1. \pl_\nu  {u_2}) {\rm dv}_{g}\\
&=-\int_{\pl M_0} (\pl_\nu v  .{f_2} - v.\pl_\nu  {u_2}){\rm dv}_{g}= -\int_{M_0} (\Delta_gv.  {u_2} - v.\Delta_g {u_2}){\rm dv}_g = 0\\
 \end{split} 
 \end{equation}
since $\Delta_g + V_2$ annihilates both $v$ and $u_2$. We substitute in the full expansion for  $u_1$ and $u_2$ and, setting $V :=V_1 - V_2$, and using 
the estimates in Lemmas \ref{termeta}, \ref{completecgo} and \ref{a0 and a1} and, we obtain 
\begin{equation}
\label{expansion}
0 = I_1 + I_2 + o(h),
\end{equation}
where 
\begin{equation}
\label{I1}
I_1 = \int_{M_0}V(a^2 + \overline{a}^2) {\rm dv}_g + 2{\rm Re}\int_{M_0}e^{2i\psi/h}V|a|^2 {\rm dv}_g,
\end{equation}
\begin{equation}
\label{I2}
 I_2 = 2h\,{\rm Re}\int_{M_0}aV\Big(e^{2i\psi/h}(\bbar{\frac{s_1}{h}+b_0}) +e^{-2i\psi/h}(\bbar{a_0+\frac{r_1}{h}}) + b_0 +  a_0+\frac{s_1+r_1}{h}\Big){\rm dv}_g.
\end{equation}

\begin{remark}\label{rem1} We observe from the last identity in Lemma \ref{termeta} that $r_1/h$ in the expression $I_2$ can be replaced by the term $\til{r}_{12}$ 
satisfying $2i\til{r}_{12}\pl \psi =b$ up to an error which can go in the $o(h)$ in \eqref{expansion}, and similarly for the term $s_1/h$ which can be replaced 
by a term $\til{s}_{12}$ independent of $h$. 
\end{remark}

We will apply the stationary phase to these two terms in the following two Lemmas.
\begin{lemma}
\label{I2 stationary phase}
The following estimates holds true 
\[I_2 = 2h\,{\rm Re}\Big(\int_{M_0}Va( b_0 +  a_0+ r_{12} + s_{12}+\til{s}_{12}+\til{r}_{12}){\rm dv}_g\Big) + o(h).\]
where $r_{12}$, $s_{12}$, $\til{r}_{12}$ and $\til{s}_{12}$ are independent of $h$.
\end{lemma}
\noindent\textbf{Proof}. We start by the following   
\begin{lemma}
\label{o of 1}
Let $f\in L^1(M_0)$, then as $h\to 0$
\[\int_{M_0}e^{2i\psi/h} f {\rm{dv}}_g = o(1).\]
\end{lemma}
\noindent{\bf Proof}. Since $C^k(M_0)$ is dense in $L^1(M_0)$ for all $k\in\nn$, it suffices to prove the Lemma 
for $f\in C^k(M_0)$.  
Let $\epsilon >0$ be small, and choose cut off function $\chi$ which is identically equal to $1$ on the boundary such that
\[\int_{M_0}\chi|f| {\rm{dv}}_g \leq \epsilon.\]
Then, splitting the integral and using stationary phase for the $1-\chi$ term, we obtain   
\[\Big|\int_{M_0}e^{2i\psi/h} f {\rm{dv}}_g\Big| \leq \Big|\int_{M_0}(1-\chi)e^{2i\psi/h} f {\rm{dv}}_g\Big| + \Big|\int_{M_0}\chi e^{2i\psi/h} f {\rm{dv}}_g\Big| \leq \epsilon + O_{\eps}(h)\]
which concludes the proof by taking $h$ small enough depending on $\eps$.\qed\\
The proof of Lemma \ref{I2 stationary phase} is a direct consequence of Lemma \ref{o of 1} and Remark \ref{rem1}.  
 \qed\\

The second Lemma will be proved in the end of this section.
\begin{lemma}
\label{I1 stationary phase}
The following estimate holds true 
\[I_1 = \int_{M_0}V(a^2 + \overline{a}^2) {\rm dv}_g + h C_pV(p)|a(p)|^2 {\rm Re}(e^{2i\psi(p)/h}) + o(h)\]
with $C_p\neq 0$ and independent of $h$.
\end{lemma}
\noindent With these two Lemmas, we can write (\ref{expansion}) as
\[ 0 = \int_{M_0}V(a^2 + \overline{a}^2) {\rm dv}_g + O(h)\]
and thus we can conclude that 
\[ 0 =  \int_{M_0}V(a^2 + \overline{a}^2) {\rm dv}_g.\]
Therefore, (\ref{expansion}) becomes
\[0 =C_pV(p)|a(p)|^2 {\rm Re}(e^{2i\psi(p)/h}) +  2{\rm Re}\Big(\int_{M_0}Va( b_0 +  s_1 +  a_0 +  r_1){\rm dv}_g\Big) + o(1).\]
Since $\psi(p)\neq 0$ we may choose a sequence of $h_j \to 0$ such that ${\rm Re}(e^{2i\psi(p)/h_j}) = 1$
and another sequence $\tilde h_j \to 0$ such that ${\rm Re}(e^{2i\psi(p)/\tilde h_j}) = -1$ for all $j$. Adding the expansion with $h = h_j$ and $h= \tilde h_j$, we deduce that 
\[ 0  = 2C_pV(p)|a(p)|^2 + o(1)\]
as $j\to \infty$, and since $C_p\neq 0$, $a(p)\neq 0$, we conclude that $V(p) = 0$. 
The set of $p\in M_0$ for which we can conclude this is dense in $M_0$ by Proposition \ref{criticalpoints}. Therefore we can conclude that $V(p) = 0$ for all $p\in M_0$.\qed\\

We now prove Lemma \ref{I1 stationary phase}. 

\noindent{\bf Proof of Lemma \ref{I1 stationary phase}}. Let $\chi$ be a smooth cutoff function on $M_0$ which is identically $1$ everywhere 
except outside a small ball containing $p$ and no other critical point of $\Phi$, and $\chi=0$ near $p$.
We split the oscillatory integral in two parts:
\[\begin{gathered}
\int_{M_0} (e^{2i\psi/h} + e^{-2i\psi/h})V |a|^2 {\rm{dv}}_g = \int_{M_0}\chi (e^{2i\psi/h} + 
 e^{-2i\psi/h})V |a|^2 {\rm{dv}}_g\\+ \int_{M_0}(1-\chi) (e^{2i\psi/h} + e^{-2i\psi/h})V |a|^2 {\rm{dv}}_g
 \end{gathered}\]
The phase $\psi$ has nondegenerate critical points, therefore, a standard application of the stationary phase at $p$ gives
\[ \int_{M_0}(1-\chi) (e^{2i\psi/h} + e^{-2i\psi/h})V(p) |a|^2 {\rm{dv}}_g = hC_p|a(p)|^2V(p){\rm Re}(e^{2i\psi(p)/h}) + o(h)\]
where $C_p$ is a non-zero number which depends on the Hessian of $\psi$ at the point $p$. 
Define the potential $\til{V}(\cdot):=V(\cdot)-V(p)\in C^{1,\alpha}(M_0)$, then we show that 
\begin{equation}\label{o of h}
\int_{M_0}(1-\chi) (e^{2i\psi/h} + e^{-2i\psi/h})\til{V}|a|^2 {\rm{dv}}_g=o(h).
\end{equation}
Indeed, first by integration by parts and using $\Delta_g\psi=0$ one has
\[\begin{split}
\int_{M_0}(1-\chi) (e^{2i\psi/h} + e^{-2i\psi/h})\til{V} |a|^2 {\rm{dv}}_g =& 
\frac{h}{2i}\int_{M_0}\cjg d(e^{2i\psi/h} - e^{-2i\psi/h}),d\psi\cjd \til{V} \frac{(1-\chi)|a|^2}{|d \psi|^2} {\rm{dv}}_g\\
=& \frac{h}{2i}\int_{M_0}(e^{2i\psi/h} - e^{-2i\psi/h})\cjg d\Big(\frac{(1-\chi) |a|^2\til{V}}{|d\psi|^2}\Big),d\psi\cjd{\rm{dv}}_g
\end{split}\]
but we can see that $\cjg d((1-\chi) |a|^2\til{V}/|d\psi|^2),d\psi\cjd\in L^1(M_0)$: 
this follows directly from the fact that $\til{V}$ is in the H\"older space $C^{1,\alpha}(M_0)$ and $\til{V}(p)=0$,
and from the non degeneracy of ${\rm Hess}(\psi)$.
It then suffice to use Lemma \ref{o of 1} to conclude that \eqref{o of h} holds.
Using similar argument, we now show that
\[\int_{M_0}\chi (e^{2i\psi/h} + e^{-2i\psi/h})V |a|^2 {\rm{dv}}_g  = o(h).\]
Indeed, since $a$ vanishes to large order at all boundary critical points of $\psi$, we may write
\[\begin{split}
\int_{M_0}\chi (e^{2i\psi/h} + e^{-2i\psi/h})V |a|^2 {\rm{dv}}_g =& \frac{h}{2i}\int_{M_0}\cjg d(e^{2i\psi/h} - e^{-2i\psi/h}),d\psi\cjd V \frac{\chi|a|^2}{|d \psi|^2} {\rm{dv}}_g\\
=&-\frac{h}{2i}\int_{M_0}(e^{2i\psi/h} - e^{-2i\psi/h}) {\rm div}_g\Big(V\frac{\chi|a|^2}{|d\psi|^2}\nabla^g\psi\Big){\rm{dv}}_g\\
&+ \frac{h}{2i}\int_{\partial M_0} (e^{2i\psi/h} - e^{-2i\psi/h})V \frac{|a|^2}{|d\psi|^2}\pl_\nu\psi\,{\rm dv}_{g}.
\end{split}\]
For the interior integral we use Lemma \ref{o of 1} to conclude that 
\[-\frac{h}{2i}\int_{M_0}(e^{2i\psi/h} - e^{-2i\psi/h}) {\rm div}_g\Big(V\frac{\chi|a|^2}{|d\psi|^2}\nabla^g\psi\Big){\rm{dv}}_g = o(h)\]
and for the boundary integral, 
we write $\partial M_0 = \Gamma_0\cup \Gamma$ and observe that on $\Gamma_0$, $\psi = 0$ so 
$(e^{2i\psi/h} - e^{-2i\psi/h}) = 0$, while on $\Gamma$ we have $V=0$ from the boundary determinacy proved in Proposition \ref{boundary determination} of the Appendix.  Therefore
\[\int_{M_0}\chi (e^{2i\psi/h} + e^{-2i\psi/h})V |a|^2 {\rm{dv}}_g = o(h)\] 
and the proof is complete.\qed

\end{section}

\section{Inverse scattering}\label{invscat}

We first obtain, as a trivial consequence of Theorem \ref{identif}, a result about inverse scattering for 
asymptotically hyperbolic surface (AH in short). Recall that an AH surface is an open complete Riemannian surface $(X,g)$ such that $X$ is the interior of a smooth compact surface with boundary $\bar{X}$, and for any smooth boundary defining function $x$ of $\pl\bar{X}$, $\bar{g}:=x^2g$ extends as a smooth metric to $\bar{X}$, with curvature tending to $-1$ at $\pl\bar{X}$.
If $V\in C^{\infty}(\bar{X})$ and $V=O(x^2)$, then we can define a scattering map as follows (see for instance \cite{JSB,GRZ} or \cite{GG}):
first the $L^2$ kernel $\ker_{L^2}(\Delta_g+V)$ is a finite dimensional subspace of $xC^{\infty}(\bar{X})$ and 
in one-to-one correspondence with $E:=\{(\pl_x\psi)|_{\pl\bar{X}}; \psi\in\ker_{L^2}(\Delta_g+V)\}$ where $\pl_x:=\nabla^{\bar{g}} x$ is the normal vector field to $\pl\bar{X}$ for $\bar{g}$, then  
for $f\in C^{\infty}(\pl\bar{X})$, there exists a function $u\in C^{\infty}(\bar{X})$, unique modulo $\ker_{L^2}(\Delta_g+V)$,
such that $(\Delta_g+V)u=0$ and $u|_{\pl\bar{X}}=f$. Then one can see that the scattering map 
$\mc{S}:C^\infty(\pl\bar{X})\to C^{\infty}(\pl\bar{X})/E$ is defined by $\mc{S}f:=\pl_x u|_{\pl\bar{X}}$. We thus obtain 
\begin{coro}\label{coroAH}
Let $(X,g)$ be an asymptotically hyperbolic manifold and let $V_1,V_2\in x^{2}C^{\infty}(\bar{X})$ be two potentials and $\Gamma\subset
\pl\bar{X}$ an open subset of the conformal boundary.  Assume that 
\[\{\pl_xu|_{\pl\bar{X}}; u\in \ker_{L^2}(\Delta_g+V_1)\}=\{\pl_xu|_{\pl\bar{X}}; u\in \ker_{L^2}(\Delta_g+V_2)\}\] 
and let $\mc{S}_j$ be the scattering map for the operator $\Delta_g+V_j$ for $j=1,2$. If 
$\mc{S}_1f=\mc{S}_2f$ on $\Gamma$ for all $f\in C_0^\infty(\Gamma)$ then $V_1=V_2$. 
\end{coro}
\noindent{\bf Proof}. Let $x$ be a smooth boundary defining function of $\pl\bar{X}$, and let $\bar{g}=x^2g$
be the compactified metric and define $\bar{V}_j:=V_j/x^2 \in C^{\infty}(\bar{X})$. 
By conformal invariance of the Laplacian in dimension $2$, one has 
\[\Delta_g+V_j=x^2(\Delta_{\bar{g}}+\bar{V}_j)\]
and so if $\ker_{L^2}(\Delta_g+V_1)=\ker_{L^2}(\Delta_g+V_2)$ and $\mc{S}_1=\mc{S}_2$ on $\Gamma$, then 
the Cauchy data spaces $\mc{C}^\Gamma_i$ for the operator $\Delta_{\bar g}+\bar{V}_j$ are the same. 
Then it suffices to apply the result in Theorem \ref{identif}.
\qed\\

Next we consider the asymptotically Euclidean scattering at $0$ frequency. 
An asymptotically Euclidean surface is a non-compact Riemann surface $(X,g)$, which compactifies into $\bar{X}$ and such that 
the metric in a collar $(0,\eps)_x\x \pl\bar{X}$ near the boundary is of the form  
\[g=\frac{dx^2}{x^4}+\frac{h(x)}{x^2}\]
where $h(x)$ is a smooth one-parameter family of metrics on $\pl\bar{X}$ with $h(0)=d\theta^2_{S^1}$ is the metric 
with length $2\pi$ on each copy of $S^1$ that forms the connected components of $\pl\bar{X}$. 
Notice that using the coordinates $r:=1/x$, $g$ is asymptotic to $dr^2+r^2d\theta^2_{S^1}$ near $r\to \infty$.
A particular case is given by the surfaces with Euclidean ends, 
ie. ends isometric to $\rr^2\setminus B(0,R)$ where $B(0,R)=\{z\in\rr^2; |z|\geq R\}$.
Note that $g$ is conformal to an asymptotically cylindrical metric, or 'b-metric' in the sense of Melrose \cite{APS}, 
\[g_b:=x^2g=\frac{dx^2}{x^2}+h(x)\]
and the Laplacian satisfies $\Delta_g=x^2\Delta_{g_b}$. Each end of $X$ is of the form $(0,\eps)_x\x S^1_{\theta}$ 
and the operator $\Delta_{g_b}$ has the expression in the ends
\[\Delta_{g_b}=-(x\pl_x)^2+\Delta_{\pl\bar{X}}+ xP(x,\theta;x\pl_x,\pl_\theta)\]
for some smooth differential operator $P(x,\theta;x\partial_x, \partial_\theta)$ in the vector fields $x\pl_x, \pl_\theta$  
down to $x=0$. 
Let us define $V_b:=x^{-2}V$, which is compactly supported and 
\[H^{2m}_b:=\{u\in L^2(X,{\rm dvol}_{g_b}); \Delta^m_{g_b}u\in L^2(X,{\rm dvol}_{g_b})\}, \quad m\in\nn_0.\] 
We also define the following spaces for $\alpha\in\rr$ 
\[F_\alpha:=\ker_{x^\alpha H^2_b}(\Delta_{g_b}+V_b).\]
Since the eigenvalues of $\Delta_{S^1}$ are $\{j^2; j\in\nn_0\}$, 
the relative Index theorem of Melrose \cite[Section 6.2]{APS} shows that $\Delta_{g_b}+V_b$ 
is Fredholm from $x^{\alpha}H^2_b$ to $x^{\alpha}H^0_b$
if $\alpha\notin \zz$.
Moreover, from subsection 2.2.4 of \cite{APS}, we have that any solution of $(\Delta_{g_b}+V_b)u=0$
in $x^{\alpha}H^2_b$ has an asymptotic expansion of the form 
\[u\sim \sum_{j>\alpha, j\in \zz}\sum_{\ell=0}^{\ell_j}x^{j}(\log x)^\ell u_{j,\ell}(\theta), \quad \textrm{ as }x\to 0\]
for some sequence $(\ell_j)_j$ of non negative integers and some smooth function $u_{j,\ell}$ on $S^1$. 
In particular, it is easy to check that $\ker_{L^2(X,{\rm dvol}_g)}(\Delta_g+V)=F_{1+\eps}$ for $\eps\in (0,1)$.
\begin{theorem} 
Let $(X,g)$ be an asymptotically Euclidean surface and $V_1,V_2$ be two compactly supported smooth potentials
and $x$ be a boundary defining function. Let $\eps\in(0,1)$
and assume that for any $j\in\zz$ and any function $\psi\in \ker_{x^{j-\eps}H^2_b}(\Delta_{g}+V_1)$ there is a 
$\varphi\in\ker_{x^{j-\eps}H^{2}_b}(\Delta_g+V_2)$ such that $\psi-\varphi=O(x^\infty)$, and conversely. 
Then $V_1=V_2$. 
\end{theorem}
\noindent\textbf{Proof}. The idea is to reduce the problem to the compact case. 
First we notice that by unique continuation, $\psi=\varphi$ where $V_1=V_2=0$.
Now it remains to prove that, if $R_\eta$ denote the restriction of smooth functions on $X$ to 
$\{x\geq \eta\}$ and $V$ is a smooth compactly supported potential in $\{x\geq \eta\}$,
then the set $\cup_{j=0}^\infty R_{\eta}(F_{-j-\eps})$ is dense in the set $N_V$ of $H^2(\{x\geq \eta\})$ solutions
of $(\Delta_g+V)u=0$. The proof is well known for positive frequency scattering (see for instance Lemma 3.2 in \cite{Mel}), 
here it is very similar so we do not give much details. The main argument is to show that it converges in $L^2$ sense
and then use elliptic regularity; the $L^2$ convergence can be shows as follows: let $f\in N_V$ such that
\[\int_{x\geq \eta}f\psi{\rm dvol}_g=0, \, \forall \, \psi \in \cup_{j=0}^\infty F_{-j-\eps},\]
then we want to show that $f=0$. By Proposition 5.64 in \cite{APS}, there exists $k\in\nn$ and a generalized right inverse $G_b$ for
$P_b=\Delta_{g_b}+V_b$ (here, as before, $x^2V_b=V$) in $x^{-k-\eps}H^2_b$, such that
$P_bG_b={\rm Id}$. This holds in $x^{-k-\eps}H^2_b$ for $k$ large enough since the cokernel of $P_b$ on this space becomes $0$ 
for $k$ large. Let $\omega=G_bf$ so that $(\Delta_{g_b}+V_b)\omega=f$, and in particular this function is $0$ 
in $\{x<\eta\}$. The asymptotic behaviour of the integral kernel $G_b(z,z')$ of $G_b$ as $z\to \infty$ 
is given in Proposition 5.64 of \cite{APS} uniformly in $z'\in\{x\geq \eta\}$, 
we have for all $J\in\nn$ and using the radial coordinates $(x,\theta)\in (0,\eps)\x S^1$ 
for $z$ in the ends
\[G_b(z,z')= \sum_{j=-k}^{J}\sum_{\ell=0}^{\ell_j} x^{j}(\log x)^{\ell}\psi_{j}(\theta,z')+o(x^{J})\]  
for some functions $\psi_{j,\ell}\in x^{k-j-\eps}H^2_b$ and some sequence $(\ell_j)_j$ of non-negative integers. 
But the fact that $(\Delta_{g_b}+V_b)G_b(z,z')=\delta(z-z')$ as distributions implies directly 
that $(\Delta_{g_b}+V_b)\psi_j(\theta,.)=0$.
Using our assumption on $f$, we deduce that $\int_{X}\psi_j(\theta,z')f(z'){\rm dvol}_{g_b}=0$ for all $j\in \nn_0$ 
and so the function $\omega$ vanish faster than all power of $x$ at infinity. Then by unique continuation, 
we deduce that $\omega=0$ in $\{x\leq \eps\}$. Since now $\omega\in H^2$, its Cauchy 
data at $x=\eta$ are $0$ and $\Delta_{g_b}+V_b$ is self adjoint for the measure ${\rm dvol}_{g_b}$, 
we can use the Green formula to obtain
\[\int_{x\geq \eta}|f|^2{\rm dvol}_{g_b}=\int_{x\geq \eta}\omega(\Delta_{g_b}+V_b)\bar{f}{\rm dvol}_{g_b}=0.\]
The $H^2$ density is easy using elliptic regularity.\qed

\begin{section}{Appendix : boundary determination}
In this appendix, we give a short proof of the fact that the partial Cauchy data space on $\Gamma\subset \pl M$ determines the potential 
on $\Gamma$ when the potential is in $C^{0,\alpha}(M)$ for some $\alpha\in(0,1)$. This result is shown for the conductivity problem on a domain of $\rr^n$ in \cite{KV}
under the assumption that the conductivity has roughly $n/2$-derivatives, it is also shown in \cite{SyU} for continuous potentials on a smooth domain of $\rr^n$ by using a limiting argument from the smooth case (which they analyze using micolocal analysis near the boundary). Alessandrini \cite{Al} also proved such a result under the assumption that the domain is Lipschitz and the coefficients of the operator are in 
$W^{1,p}$ for $p>n$,  while Brown \cite{Br}Ê studied the case of Lipschitz domains with a continuous conductivity.   
Since the result in our setting is not explicitly written down, but certainly known from specialists, we provide a short proof 
without too many details, and using the approach of \cite{Br}. 
We shall prove 
\begin{proposition} \label{boundary determination}
Let $\Gamma\subset \pl M_0$ be an non-empty open subset of the boundary. If $V_1,V_2\in C^{0,\alpha}(M)$ for some $\alpha>0$ and their associated Cauchy data spaces ${\cal C}^\Gamma_1,{\mc C}^\Gamma_2$ defined in \eqref{cauchydata} are equal, then $V_1|_{\Gamma}=V_2|_{\Gamma}$.
\end{proposition}
The key to proving this proposition is the existence of solutions to $(\Delta_g + V_i) u = 0$ which concentrate near a point $p\in \Gamma$. First we need a solvability result for the equation $(\Delta_g+V_i)u=f$, which is an easy consequence of the Carleman estimate of Proposition \ref{carlemanestimate}, and follows the method of Salo-Tzou \cite[Section 6]{SaT}.  
If we fix $h>0$ small and take $\varphi=1$ in the Carleman estimate of Proposition \ref{carlemanestimate}, we obtain easily that there is a constant $C$ such that for all functions in $H^2(M_0)$ satisfying $u|_{\pl M_0}=0$
\begin{equation}
\label{ellipticforschrodinger}
\|u\|^2_{H^2} + \|\pl_\nu u\|^2_{L^2(\Gamma_0)}\leq C (\|(\Delta +V_i) u\|^2_{L^2}+ \|\pl_\nu u\|^2_{L^2(\Gamma)}) 
\end{equation}
As a consequence, we deduce the following solvability result: let 
$${\cal B} := \{w \in H^2(M_0) \cap H^1_0(M_0) \mid \pl_\nu w|_{\Gamma} = 0\}$$ 
be the closed subspace of $H^2(M_0)$ under the $H^2$ norm and let ${\cal B}^*$ be its dual space then
\begin{coro}
\label{solvability}
Let $i=1,2$, then for all $f\in L^2(M_0)$ there exists $u\in H^2(M_0)$ solving the equation 
\[(\Delta_g + V_i) u = f\]
with boundary condition $u|_{\Gamma_0} = 0$, and $\|u\|_{L^2}\leq C\|f\|_{{\cal B}^*}$.
\end{coro}
\noindent{\bf Proof}.
Set ${\cal A} := \{ w\in H^1_0(M_0)\mid (\Delta_g + V_i) w \in L^2, \pl_\nu w|_{\Gamma} = 0\}$ equipped with the inner product
\[(v,w)_{\cal A} := \int_{M_0}(\Delta_g + V_i)v(\Delta_g +V_i)\bar w \, {\rm dv}_g.\]
Thanks to \eqref{ellipticforschrodinger}, $\cal A$ is a Hilbert space and $\cal A = \cal B$. For each $f\in L^2(M_0)$, 
let us define the linear functional on $\cal B$
\[L_f : w \mapsto \int_{M_0} wf \, {\rm dv}_g.\]
By  \eqref{ellipticforschrodinger}, we have that for all $w\in \cal B$
\[|L_f(w)| \leq \|f\|_{{\cal B}^*}\|w\|_{\cal B} \leq  \|f\|_{{\cal B}^*}\|w\|_{\cal A}.\]
Therefore, by Riesz Theorem, there exists $v_f\in \cal A$ such that 
\[\int_{M_0} (\Delta_g + V) \bbar {v_f} (\Delta_g + V_i)w\, {\rm dv}_g = \int_{M_0} wf {\rm dv}_g\]
for all $w\in\cal B$. Furthermore, $\|(\Delta_g +V_i) v_f\|_{L^2} \leq \|f\|_{\cal B^*}$. Setting $u := (\Delta_g + V_i) \bbar {v_f}$ we have that $(\Delta_g + V_i) u = f$ and $\|u\|_{L^2} \leq \|f\|_{\cal B^*}$. To obtain the boundary condition for $u$, observe that since
\[\int_{M_0}u(\Delta_g + V_i) w dv_g = \int_{M_0} f w \, {\rm dv}_g\]
for all $w\in\cal B$, by Green's theorem
\[\int_{M_0}u \pl_\nu w\, {\rm dv}_g = 0 = \int_{\Gamma_0} u\pl_\nu w \, {\rm dv}_g\]
for all $w\in\cal B$. This implies $u = 0$ on $\Gamma_0$. 
\qed\\

Clearly, it suffices to assume that $\Gamma$ is a small piece of the boundary which is contained in a single coordinate chart with complex coordinates $z=x+iy$ where $|z|\leq 1$,  ${\rm Im}(z)>0$ and the boundary is given by $\{y=0\}$. Moreover the metric is of the form
$e^{2\rho}|dz|^2$ for some smooth function $\rho$. Let $p\in \Gamma$ and possibly by translating the coordinates, we can assume that $p=\{z=0\}$. Let $\eta \in C^\infty(M_0)$ be a cutoff function supported in a small neighbourhood of $p$. 
For $h>0$ small, we define the the smooth function $v_h\in C^\infty(M_0)$ supported near $p$ via the coordinate chart $Z=(x,y)\in\rr^2$ by 
\begin{equation}\label{vh}
v_h(Z):= \eta(Z/\sqrt{h})e^{\frac{1}{h}\alpha\cdot Z}
\end{equation} 
where $\alpha:=(i,-1) \in\cc^2$ is chosen such that $\alpha\cdot\alpha = 0$. We thus get $(\pl_{x}^2 + \pl_{y}^2)e^{\alpha\cdot Z} = 0$ and thus $\Delta_ge^{\alpha\cdot Z}=0$ by conformal covariance of the Laplacian. 
Therefore, we have in local coordinates
\begin{equation}\label{error term}
\Delta_g v_h(Z) = \frac{1}{h}e^{\frac{1}{h}\alpha\cdot Z}(\Delta_g\eta)(\frac{Z}{\sqrt h})  + \frac{2}{h^{3/2}}e^{\frac{1}{h}\alpha\cdot Z}\cjg d\eta(\frac{Z}{\sqrt h}),  \alpha.dZ\cjd_g.
\end{equation}
\begin{lemma}
\label{solutions concentrating near boundary}
If $V\in C^{0,\alpha}(M_0)$ for $q>2$, then there exists a solution $u_h\in H^2$ to $(\Delta_g + V)u = 0$ of the form
\[u_h = v_h + R_h,\]
with $v_h$ defined in \eqref{vh} and $\|R_h\|_{L^2} \leq C h^{5/4}$, satisfying ${\rm supp}(R_h|_{\pl M_0})\subset\Gamma$.
\end{lemma}
\noindent {\bf Proof of Lemma \ref{solutions concentrating near boundary}}. 
We need to find $R_h$ satisfying $\|R_h\|_{L^2} \leq C h^{5/4}$ and solving 
\[(\Delta_g + V) R_h = -(\Delta_g + V)v_h =: M_h.\] 
Thanks to Corollary \ref{solvability}, it suffices to show that $\|M_h\|_{\cal B^*}\leq Ch^{5/4}$.  Thus, let $w\in \cal B$, then we have by (\ref{error term})
\[\int_{M_0} w M_h\, {\rm dv}_g =  I_1 + I_2 + I_3\]
where 
\[I_1 := \int_{|Z|\leq\sqrt{h}} wV_i  v_h e^{2\rho}dZ , \quad I_2 := \frac{1}{h} \int_{|Z|\leq \sqrt{h}}w e^{\frac{1}{h}\alpha\cdot Z}\chi_1(\frac{Z}{\sqrt h} )e^{2\rho}dZ ,\]
\[I_3 := \frac{2}{ h^{3/2}} \int_{|Z|\leq \sqrt{h}} w \chi_2(\frac{Z}{\sqrt h}) e^{\frac{1}{h}\alpha\cdot Z}\, dZ\]
and $\chi_1 = \Delta_g \eta$, $\chi_2 = i\pl_x\eta - \pl_y\eta$.
In the above equation the third term $I_3$ has the worst growth when $h\to 0$. We will analyze its behavior and the preceding terms can be treated in similar fashion. One has
\[\begin{split} 
I_3 =& -h^{1/2}\int_{|Z|\leq \sqrt{h}}w \chi_2(\frac{Z}{\sqrt h})(\pl_{x}- i \pl_{y})^2 e^{\frac{1}{h}\alpha\cdot Z}dZ \\ 
=&- h^{1/2}\int_{|Z|\leq \sqrt{h}} (\pl_{x}- i \pl_{y})^2 (w \chi_2(\frac{Z}{\sqrt h} ))e^{\frac{1}{h}\alpha\cdot x} dZ.
\end{split}\]
Notice that the boundary term in the integration by parts vanishes because $w\in H^1_0$ and $\pl_\nu w|_{\pl M_0}$ vanishes on the support of $\eta$. The term $(\pl_{x}- i \pl_{y})^2 (w \chi_2(\frac{Z}{\sqrt h}))$ has derivatives hitting both $\chi_2(\frac{Z}{\sqrt h})$ and $w$. The worst growth in $h$ would occur when both derivatives hit $\chi_2(\frac{Z}{\sqrt h})$ in which case a $h^{-1}$ factor would come out. Combined with the $h^{1/2}$ term in front of the integral this gives a total of a $h^{-1/2}$ in front. By this observation we have improved the growth from $h^{-3/2}$ to $h^{-1/2}$. Repeating this line of argument and using Cauchy-Schwarz inequality, we can see that 
$|I_3| \leq C h^{5/4}\|w\|_{H^2}$ (an elementary computation shows that functions of the form $\chi(Z/\sqrt{h})e^{\frac{1}{h}\alpha\cdot Z}$ have $L^2$ norm bounded by $Ch^{3/4}$).
Therefore, $\|\frac{1}{h^{3/2}}\chi_2(\frac{Z}{\sqrt h}) e^{\frac{1}{h}\alpha\cdot Z}\|_{{\cal B}'} \leq Ch^{5/4}$ and we are done.
\qed\\

\noindent {\bf Proof of Proposition \ref{boundary determination}}. It suffices to plug the solutions $u_h^1,u_h^2$ 
from Proposition \ref{solutions concentrating near boundary} into the boundary integral identity \eqref{integralid}. 
A simple calculation using the fact that $V_1-V_2$ is in $C^{0,\gamma}(M_0)$  yields that 
\[0= \int_{M_0} u_h^1(V_1-V_2) u^2_h \,{\rm dv}_g = C h^{3/2} (V_1(p)-V_2(p)) + o(h^{3/2})\]
and we are done.\qed

\end{section}

\subsection{Acknowledgements}
This work started during a summer evening in Pisa thanks to the hospitality of M. Mazzucchelli and A.G. Lecuona. We thank Sam Lisi, Rafe Mazzeo, Mikko Salo, Eleny Ionel for pointing out very helpful references. 
C.G. thanks MSRI and the organizers of the 'Analysis on Singular spaces' 2008 program  
for support during part of this project. Part of this work was done while C.G. was 
visiting IAS under an NSF fellowship number No. DMS-0635607. 
L.T is supported by NSF Grant No. DMS-0807502.

\end{document}